\documentclass[10pt]{amsart}
\usepackage{amsfonts,amssymb}
\RequirePackage{ifpdf}
\usepackage{amsmath}
\usepackage{color}
\usepackage{pstcol,pst-node}
\usepackage{graphics}
\usepackage{epic}
\usepackage{curves}
\usepackage{multirow}

\def\appendix#1{
\addtocounter{section}{1} \setcounter{equation}{0}
\renewcommand{\thesection}{\Alph{section}}
\section*{Appendix \thesection\protect\indent\quad
#1}
}
\renewcommand{\theequation}{\thesection.\arabic{equation}}

\catcode`\@=11
\def\marginnote#1{}

\newcount\hour
\newcount\minute
\newtoks\amorpm
\hour=\time\divide\hour by60 \minute=\time{\multiply\hour by60
\global\advance\minute by-\hour}
\edef\standardtime{{\ifnum\hour<12 \global\amorpm={am}%
        \else\global\amorpm={pm}\advance\hour by-12 \fi
        \ifnum\hour=0 \hour=12 \fi
        \number\hour:\ifnum\minute<10 0\fi\number\minute\the\amorpm}}
\edef\militarytime{\number\hour:\ifnum\minute<100\fi\number\minute}

\newcommand{\mapi}{{\stackrel{{i}}{{\to}}}}

\def\draftlabel#1{{\@bsphack\if@filesw {\let\thepage\relax
      \xdef\@gtempa{\write\@auxout{\string
          \newlabel{#1}{{\@currentlabel}{\thepage}}}}}\@gtempa \if@nobreak
    \ifvmode\nobreak\fi\fi\fi\@esphack} \gdef\@eqnlabel{#1}}
    \def\@eqnlabel{}
\def\@vacuum{}
\def\draftmarginnote#1{\marginpar{\raggedright\scriptsize\tt#1}}

\def\draft{
  \oddsidemargin -.5truein
  \def\@oddfoot{\footnotesize \sl preliminary draft \hfil
    \rm\thepage\hfil\sl\today\quad\militarytime}
  \let\@evenfoot\@oddfoot \overfullrule 3pt
    \let\label=\draftlabel
    \let\marginnote=\draftmarginnote
  \def\@eqnnum{(\theequation)\rlap{\kern\marginparsep\tt\@eqnlabel}%
    \global\let\@eqnlabel\@vacuum}

  }

\def\be{\begin{equation}}
\def\ee{\end{equation}}
\def\bea{\begin{eqnarray}}
\def\eea{\end{eqnarray}}
\def\<{\langle}
\def\>{\rangle}
\def\nn{\nonumber}

\def\otim{\mathop{\otimes}}

\def\ocomma{{\phantom{\Bigm|}^{\phantom {X}}_{\raise-1.5pt\hbox{,}}\!\!\!\!\!\!\otimes}}

\newcommand{\sheet}[2]{{\stackrel{{#1}}{{#2}}}}
\newcommand{\groupoid}{{\stackrel{{\rightarrow}}{{_\rightarrow}}}}

\newtheorem{theorem}{Theorem}[section]
\newtheorem{lm}[theorem]{Lemma}
\newtheorem{prop}[theorem]{Proposition}
\theoremstyle{definition}
\newtheorem{defin}[theorem]{Notation}

\newtheorem{remark}[theorem]{Remark}

\newtheorem{conjecture}[theorem]{Conjecture}

\allowdisplaybreaks

\begin{document}

\title[Poisson  space of bilinear forms]
{On a Poisson  space of bilinear forms with a Poisson Lie action}
\author{Leonid Chekhov$^{\ast}$}\thanks{$^{\ast}$Steklov Mathematical Institute and  Laboratoire Poncelet,
Moscow, Russia; School of Mathematica, Loughborough University, UK. Email: chekhov@mi.ras.ru.}
\author{Marta Mazzocco$^\dagger$}\thanks{$^\dagger$School of Mathematica, Loughborough University, LE11 3TU UK. m.mazzocco@lboro.ac.uk}

\begin{abstract}
We consider the space $\mathcal A$ of bilinear forms on $\mathbb C^N$
with defining matrix $\mathbb A$  endowed with the quadratic Poisson structure studied  by the authors in \cite{ChM}.
We classify all possible quadratic brackets
on $(B,\mathbb A)\in GL_N\times \mathcal A$ with the property that the natural action $\mathbb A\mapsto B\mathbb AB^{\text{T}}$ of the $GL_N$ Poisson--Lie group on the space $\mathcal A$ is a Poisson action thus endowing $\mathcal A$
with the structure of Poisson space. Beside the product Poisson structure on $GL_N\times \mathcal A$ we find two more
(dual to each other) structures for which (in contrast to the product Poisson structure) we can implement
the reduction to the space of bilinear forms with block upper triangular defining  matrices by Dirac procedure. We consider the generalisation of the above construction to triples $(B,C,\mathbb A)\in GL_N\times GL_N\times \mathcal A$
with the Poisson action $\mathbb A\mapsto B\mathbb AC^{\text{T}}$ and show that $\mathcal A$ then acquires the structure of  Poisson symmetric space. We study also the generalisation  to chains of transformations and to the
quantum and quantum affine algebras {and the relation between the construction of Poisson symmetric spaces
and that of the Poisson groupoid.}
\end{abstract}

\maketitle

\section{Introduction}
In this paper, we identify bilinear forms on $\mathbb C^N$
$$
\langle x,y\rangle:= x^\text{T} \mathbb A y,\qquad \forall\, x,\,y\in\mathbb C^N,
\qquad \mathbb A\in \hbox{Mat}_N(\mathbb C),
$$
with their defining matrix $\mathbb A$. We denote by $\mathcal A$ the space of such bilinear forms.

In~\cite{ChM}, the authors studied a quadratic Poisson algebra structure on the
space $\mathcal A$ of bilinear forms on $\mathbb C^{N}$ with the property that
for any $n,m\in\mathbb N$ such that $n m =N$, the restriction of the
Poisson algebra to the space $\mathcal A_{n,m}$ of bilinear forms with block-upper-triangular (b.u.t.) defining
matrix composed from blocks of size $m\times m$ is Poisson:
\bea
\label{Poisson}
\{a_{i,j},a_{k,l}\}&=&\bigl({\rm sign}(j-l)+{\rm sign}(i-k)\bigr)a_{i,l}a_{k,j}+\\
&&
+\bigl({\rm sign}(j-k)+1\bigr)a_{j,l}a_{i,k}
+\bigl({\rm sign}(i-l)-1\bigr)a_{l,j}a_{k,i}.\nn
\eea
These algebras were studied previously in the upper-triangular case in \cite{GK}, \cite{NR}, \cite{NRZ}
and in the case of $2\times2$ blocks in \cite{MR}, \cite{MRS} (see also monograph \cite{Molev}) in relation to
various algebraic and geometric systems.

In the $r$-matrix notation explained in Appendix~A the same bracket (\ref{Poisson})
can be written in a more concise form:
\be
\label{Poisson-r}
\{\sheet1{\mathbb A}\ocomma\sheet2{\mathbb A}\}=r_{12}(\sheet1{\mathbb A}\otimes\sheet2{\mathbb A})
-(\sheet1{\mathbb A}\otimes\sheet2{\mathbb A})r_{12}+\sheet1{\mathbb A}r_{12}^{t_1}\sheet2{\mathbb A}
-\sheet2{\mathbb A}r_{12}^{t_1}\sheet1{\mathbb A},
\ee
where
\be
r_{12}=\sum_i \sheet1{E}_{i,i}\otimes \sheet2{E}_{i,i}+2\sum_{i>j} \sheet1{E}_{i,j}\otimes \sheet2{E}_{j,i}
=2\sum_{i,j}\theta(i-j)\sheet1{E}_{i,j}\otimes \sheet2{E}_{j,i},
\label{r-matrix}
\ee
with $\theta(x)=\{1,x>0,\, 1/2,x=0,\,0,x<0\}$
is the classical (trigonometrical) $r$-matrix.

It is natural to consider the following action of $GL_N$ on the space of bilinear forms $\mathcal A$:
\be
\forall B\in GL_N,\quad  B: {\mathbb A}\mapsto {\mathbb A}':= B{\mathbb A}{B^{\text{T}}}.
\label{BAB}
\ee
In \cite{ChM}, the groupoid $\Gamma$ of morphisms of the space  $\mathcal A_{n,m}$ of  b.u.t.
bilinear forms was defined in such a way that all morphisms automatically preserve the Poisson algebra on $\mathcal A_{n,m}$. {The condition for transformation (\ref{BAB}) to be Poisson is in a sense opposite to the
standard construction of  a Poisson (symplectic) groupoid \cite{Karasev},~\cite{Weinstein} in which
the source $s:({\mathbb A},B)\to {\mathbb A}$ and target $t:({\mathbb A},B)\to B{\mathbb A}B^{\text{T}}$
projections are respectively anti-Poisson and Poisson. We present our treatment of a (possibly more familiar to
the reader) Poisson groupoid construction in Sec.~\ref{s:groupoid}.}

In the main part of this paper we show that it is in fact more natural to ask the following question along the Jiang-Hua Lu  approach \cite{Lu}: considering the action of the $GL_N$ Poisson Lie group on $\mathcal A$
we must classify all possible Poisson brackets on $GL_N\times\mathcal A$ such that this action is Poisson and $\mathcal A$ is a Poisson space in the sense of \cite{Dr}.

We assume that the bracket between $B$-matrices has the standard Lie--Poisson form
\be
\{\sheet1 B\ocomma\sheet2 B\}=r_{12} (\sheet1 B\otimes\sheet2 B) - (\sheet1 B\otimes\sheet2 B) r_{12}, \label{BB-bracket}
\ee
and prove (see Theorem \ref{lm-generalization}) that in order for the action of $GL_N$ on $\mathcal A$ to be Poisson, the brackets between entries
of the ${\mathbb A}$-matrix and entries of the $B$-matrix must  have a special quadratic form
\be
\{\sheet1{B}\ocomma\sheet2{\mathbb A}\}=\sheet1{B}Q_{12}\sheet2{\mathbb A}+\sheet1{B}\sheet2{\mathbb A}Q_{12}^{t_2},
\label{bracket-Q}
\ee
 with the matrix $Q_{12}$ taking only three possible choices:
 \begin{itemize}
\item[(i)] $Q_{12}=0$;
\item[(ii)] $Q_{12}=-r_{12}^{t_2}$;
\item[(iii)] $Q_{12}=r_{12}^{t_1}$.
	\end{itemize}
In the text, we refer to these choices by their numbers (i), (ii), and (iii). The Poisson structure (i) on $GL_N\times\mathcal A$  is the usual product Poisson structure (i.e., with $Q_{12}=0$), while the other two are dual to each other (see Lemma \ref{duality-ii-iii}). 

Stress that in none of the above cases $GL_N\times\mathcal A\,\groupoid\,\mathcal A$  is a Poisson groupoid (in particular the target map $\beta:(B,\mathbb A) \to  B{\mathbb A}{B^{\text{T}}}$ is not an anti-Poisson map). We can nevertheless
endow $GL_N\times GL_N$ with a Poisson symmetric Lie group structure and show that $\mathcal A$ is the symmetric space associated to it in the sense defined by Fernandes in \cite{F2} (see Theorem \ref{th:p-sym} here below). Loosely speaking this means that we can think of $\mathcal A$ as the set of matrices $\mathbb A$ of the form $\mathbb A= BC^{\text{T}}$, where $(B,C)\in GL_N\times GL_N$, and
we prove that the Poisson bracket (\ref{Poisson-r}) between entries of $\mathbb A$ is in fact induced by the Poisson structure on  $GL_N\times GL_N$ by the identification $\mathbb A= BC^{\text{T}}$.

Again, we classify all quadratic Poisson brackets on $GL_N\times GL_N\times\mathcal A$ for which the Poisson Lie group action on $\mathcal A$  defined by
$$
(\mathbb A,B,C) \mapsto B \mathbb A C^{\text{T}}, \qquad\forall  (B,C)\in GL_N\times GL_N
$$
is Poisson (see Lemma \ref{lem-ABC-new}).

This raises the following question: if we identify $\mathbb A= BC^{\text{T}}$, is it true that the Poisson bracket on
$GL_N\times GL_N\times\mathcal A$ is induced from the bracket among entries of $B$ and $C$? It turns out that the correct way to formulate (and indeed answer) this question is in terms of {\it chains:}\/ we introduce Poisson brackets on  $\bigotimes_1^K  (GL_N\times GL_N)\otimes\mathcal A$ such that the action of the Poisson Lie group
$\bigotimes_K (GL_N\times GL_N)$ on $\mathcal A$ is Poisson for every $K$. Then, if we identify $\mathbb A= B_1C_1^{\text{T}}$, where $(B_1,C_1)$ is an element  in the first copy of $GL_N\times GL_N$,  the Poisson brackets on  $\bigotimes_2^{K}  (GL_N\times GL_N)\otimes\mathcal A$ are induced by the one on $\bigotimes_1^{K}  (GL_N\times GL_N)\otimes\mathcal A$  (see Section \ref{s:groupoid-new-BC}).

In Section \ref{s:central} we classify all central elements for all Poisson brackets on $GL_N\times\mathcal A$, $GL_N\times GL_N\times\mathcal A$, and on chains of $B$-matrices and $(B,C)$-pairs.

The next set of results deals with the natural question of reductions to the space  $\mathcal A_{n,m}$ of b.u.t. bilinear forms. In order to implement this reduction, we need to introduce the set of constraints
(\ref{b.u.t.-constraints-A}) and (\ref{b.u.t.-constraints-BABt})
on blocks $A_{I,J}$ and $A'_{I,J}$ of the respective matrices $\mathbb A$ and $\mathbb A':=B\mathbb A B^{\text{T}}$.

This is where the non trivial  Poisson brackets (ii) or (iii) between $A$ and $B$ become necessary. In fact, these constraints do not Poisson commute with all other elements on the constraint surface, therefore this reduction is not Poisson and a Dirac reduction is needed. However in order for the Dirac reduction to work, the matrix given by all Poisson brackets between constraints must be non-degenerate on the constraint surface, and this is only possible when $A$ and $B$ do not Poisson commute with each other.

We illustrate in detail how the Dirac procedure works in the case of upper--triangular matrices $\mathbb A$, i.e. on the space  $\mathcal A_{n,1}$.

It is interesting to observe that in this case we can solve the constraint equations for $\mathbb A$ thus obtaining the entries  $a_{i,j}$ as functions $F_{i,j}[B]$. In this way, we identify ${\mathbb A}$ with the new upper-triangular matrix ${\mathbb F}$. {In what follows, we can proceed in two, very different, ways. In our original formulation of the Poisson   space, we can treat equations ${\mathbb A}_{i,j}-F_{i,j}[B]=0$, $i<j$,
${\mathbb A}_{i,j}-\delta_{i,j}=0$, $i\ge j$, as an equivalent set of second-kind constraints implementing the
same Dirac procedure as above (Sec.~\ref{s:reductions}). On the other hand, we can just induce the Poisson brackets
on ${\mathbb A}$ and between $\mathbb A$ and $B$ from the Lie Poisson brackets on $B$, which, as we demonstrate in
Sec.~\ref{s:groupoid} results in the structure of Poisson (symplectic) groupoid on the pair $({\mathbb A},B)$.}

In Secs.~\ref{s:quantization} and \ref{s:quant-spectral}, we use the r-matrix formalism to quantise all brackets and produce a quantum affine version of the Poisson algebra  (ii) on $GL_N\times \mathcal A$ and prove the quantum Jacobi relations (the formulae for the Poisson algebra (iii) can be deduced by duality).

We finally address the symplectic groupoid construction by \cite{Karasev},~\cite{Weinstein}, \cite{Bondal}.
Assuming that the brackets on $F_{i,j}[B]$ are induced by those on the entries of $B$ we obtain that
$F_{i,j}[B]$ satisfy the same relations as the entries of $\mathbb A$ with opposite sign. The thus found brackets on the set of $({\mathbb F},B)$ pairs are again quadratic and admit an $r$-matrix form of writing. We can therefore
extend these brackets to the general case $({\mathbb F},B)\in GL_N\times GL_N$. We find that
the mapping ${\mathbb F}\mapsto \tilde {\mathbb F}:=B{\mathbb F}B^{\text{T}}$
is then indeed an {\it anti}automorphism of the
Poisson algebra for
${\mathbb F}$ whereas all entries of ${\mathbb F}$ and $B{\mathbb F}B^{\text{T}}$ mutually Poisson commute.
This is in accordance with the factorization property of the sympectic groupoid \cite{Weinstein},~\cite{kirill}.
We then show that this alternative bracket admits upper-triangular and block-upper-triangular Poisson reductions
without involving the Dirac procedure and looks therefore quite
attractive on the first sight. Its disadvantage, to our opinion, is that it
does not produce a nontrivial dynamics resulting
just in the appearance of two separate copies of the original Poisson algebra for ${\mathbb A}$ (with opposite signs) sharing the same central elements generated by
$\det({\mathbb F}+\lambda {\mathbb F}^{\text{T}})=\det(\tilde{\mathbb F}+\lambda \tilde{\mathbb F}^{\text{T}})$.
It nevertheless satisfy the definition of the Poisson (symplectic) groupoid \cite{Mikami-Wein}.

\begin{remark}
It si worth noting that $\mathcal A$ arises naturally as the unipotent radical of Borel subgroups of complex simple Lie groups, can be identified with Schubert cells on flag varieties. In this way the setting of our paper can be related to Goodearl and Yakimov work \cite{GY}. The investigation on the exact relationship between the Poisson structures in that paper and ours in postponed to subsequent publication.\footnote{We thank the anonymous referee for this observation}
\end{remark}

\section{The Poisson algebra  on $GL_N\times\mathcal A$}
\setcounter{equation}{0}

In this section our aim is to find a Poisson structure on $GL_N\times\mathcal A$ such that the $GL_N$-action
\be
\label{groupoid}
{\mathbb A}\mapsto {\mathbb A}':= B{\mathbb A}{B^{\text{T}}},
\ee
on the space  ${\mathcal A}$ is Poisson.

We assume that the bracket (\ref{Poisson-r}) (or (\ref{Poisson}) in the coordinate form of writing) holds
on the space of ${\mathbb A}$-matrices and look for
such brackets $\{\sheet1{\mathbb A}\ocomma\sheet2 B\}$ and $\{\sheet1 B\ocomma\sheet2 B\}$
that
\bea
\{(\sheet1{B\mathbb A {B^{\text{T}}}})\ocomma(\sheet2{B\mathbb A B^{\text{T}}})\}
&=&r_{12}(\sheet1{B\mathbb A B^{\text{T}}})\otimes(\sheet2{B\mathbb A B^{\text{T}}})
-(\sheet1{B\mathbb A B^{\text{T}}})\otimes(\sheet2{B\mathbb A B^{\text{T}}})r_{12}+\nn\\
&{}&+(\sheet1{B\mathbb A B^{\text{T}}})r_{12}^{t_1}
(\sheet2{B\mathbb A B^{\text{T}}})-(\sheet2{B\mathbb A {B^{\text{T}}}})r_{12}^{t_1}(\sheet1{B\mathbb A {B^{\text{T}}}}).
\label{bracket-transformed}
\eea
We naturally assume that all the brackets are quadratic and preserve the number of $a$ and $b$ items.
We also assume that the bracket between $B$-matrices has the standard Lie--Poisson form (\ref{BB-bracket}),
or in coordinates:
\be
\label{Poisson-b}
\{b_{i,j},b_{k,l}\}=\bigl({\rm sign}(j-l)+{\rm sign}(i-k)\bigr)b_{i,l}b_{k,j}.
\ee

\begin{theorem}\label{lm-generalization}
Given the Poisson brackets (\ref{Poisson-r}) between entries
of the ${\mathbb A}$-matrix and (\ref{Poisson-b}) between
entries of the $B$-matrix, we have exactly three choices for the quadratic brackets between
$\mathbb A$ and $B$ such that (a) the
mapping $\mathbb A\mapsto B\mathbb A{B^{\text{T}}}$
is an automorphism of the Poisson algebra (\ref{Poisson-r}) and (b) the bracket satisfies
all the Jacobi relations: all these brackets have the form (\ref{bracket-Q})
with
\begin{itemize}
\item[(i)] $Q_{12}=0$;
\item[(ii)] $Q_{12}=-r_{12}^{t_2}$;
\item[(iii)] $Q_{12}=r_{12}^{t_1}$.
\end{itemize}
\end{theorem}

\proof
We begin with the observation that if we begin with evaluating the brackets between entries of
the matrices $B$ in (\ref{bracket-transformed}) then, among eight terms four will produce the right-hand side of this relation whether all the remaining terms have to have the structure
$$
\sheet1{B}\sheet2{B}\{\dots\} \sheet1{B^{\text{T}}}\sheet2{B^{\text{T}}},
$$
where the ellipses stand for a combination of $r$-matrices and elements of matrices ${\mathbb A}$.
The result of evaluation of brackets between entries of matrices ${\mathbb A}$ in
(\ref{bracket-transformed}) has the same form: we obtain expressions sandwiched between
$\sheet1{B}\sheet2{B}$ and $\sheet1{B^{\text{T}}}\sheet2{B^{\text{T}}}$. We therefore assume that,
in order to be able to attain proper cancellations, the inter-brackets between $B$ and ${\mathbb A}$
in the expression (\ref{bracket-transformed}) must result in the same sandwiched structure. We therefore look for  brackets having a (general) quadratic form
\be
\{\sheet1{B}\ocomma\sheet2{\mathbb A}\}=\sheet1{B}Q_{12}\sheet2{\mathbb A}+\sheet1{B}\sheet2{\mathbb A}R_{12}.
\label{quadratic}
\ee
The direct calculation then shows that
in order to yield an automorphism of the Poisson algebra (\ref{Poisson-r})
the matrices $Q_{12}$ and $R_{12}$ must satisfy the conditions
$$
R_{12}=Q_{21}^{\text{T}}\quad\hbox{and}\quad Q_{12}-Q_{21}=\kappa P_{12},
$$
where $P_{12}=\sum_{i,j}E_{i,j}\otimes E_{j,i}$ is the standard permutation matrix and $\kappa$ is a constant.
Substituting this anzatz and verifying all the Jacobi relations in the triples
$(\sheet1{B},\sheet2{B},\sheet3{\mathbb A})$ and $(\sheet1{B},\sheet2{\mathbb A},\sheet3{\mathbb A})$
we observe that, first, $\kappa=0$ and therefore $Q_{12}=Q_{21}$, and, second, that $Q_{12}$ may assume only the
above three forms provided the $(B,B)$-brackets and the
$({\mathbb A},{\mathbb A})$-brackets are given by the respective
formulas (\ref{Poisson-r}) and (\ref{BB-bracket}).
\endproof

Although the brackets (ii) and (iii) look less natural that (i), they manifest interesting symmetries as shown by the following result:

\begin{prop}\label{prop-anti}
Assuming that $\mathbb A\in GL_N(\mathbb C)$ and $B \in GL_N(\mathbb C)$, we obtain that
in the cases (ii) and (iii) of brackets (\ref{bracket-Q}),
the quantities
\be
{\mathfrak A}:=B{\mathbb A}^{-\text{T}}{B^{\text{T}}}
\label{mathfrakA}
\ee
satisfy the {\bf same} Poisson algebra (\ref{Poisson-r}) as both $\mathbb A$ and $B{\mathbb A}{B^{\text{T}}}$.
\end{prop}

\proof  Straightforward calculation using the $r$-matrix form of writing for the
corresponding brackets.\endproof

In order to save the space, we are mostly dealing with the bracket (ii) in what follows; the case of
the bracket (iii) is in fact dual to it as proved in the following:

\begin{lm}\label{duality-ii-iii}
For $\mathbb A$ and $B$ from $GL_N(\mathbb C)$ we have the following (anti)homomorphism between
the $(\mathbb A,B)$-algebras (ii) and (iii). For the
quantities $\mathbb A':=\mathbb A^{-1}$ and $B'=B^{-\text{T}}$ the brackets are as follows: the bracket
$\{\sheet1{\mathbb A'}\ocomma\sheet2{\mathbb A'}\}$ has form (\ref{Poisson-r}) with the overall minus sign,
the bracket $\{\sheet1{B'}\ocomma\sheet2{B'}\}$ has form (\ref{BB-bracket}) with the overall minus sign whereas
the bracket $\{\sheet1{B'}\ocomma\sheet2{\mathbb A'}\}$ has form (\ref{bracket-Q}) with the minus sign and
with the matrix $Q'_{12}$ of type (iii) if the matrix $Q_{12}$ was of type (ii) and vice versa. (Of course,
if $\mathbb A$ and $B$ Poisson commute then do the matrices $\mathbb A'$ and $B'$.)
\end{lm}

\proof  Straightforward calculation using the $r$-matrix form of writing for the
corresponding brackets.\endproof

\section{The Poisson algebra  on $GL_N\times GL_N \times\mathcal A$}\label{ss:Lie-Poisson}
\setcounter{equation}{0}

We now consider the general transformation
\be
\mathbb A\mapsto B\mathbb A C^{\text{T}},
\label{triple}
\ee
for which we find the brackets between the matrices $B$ and $C$ that
preserve the Poisson relations (\ref{Poisson-r}).

We endow the product $GL_N\times GL_N$ with the following brackets: $\forall (B,C)\in GL_N\times GL_N$
\bea
&{}&\{\sheet1 B\ocomma\sheet2 B\}=r_{12} (\sheet1 B\otimes\sheet2 B) - (\sheet1 B\otimes\sheet2 B) r_{12}, \label{BB-bracket1}\\
&{}&\{\sheet1 C\ocomma\sheet2 C\}=r_{12} (\sheet1 C\otimes\sheet2 C) - (\sheet1 C\otimes\sheet2 C) r_{12}, \label{CC-bracket1}\\
&{}&\{\sheet1 C\ocomma\sheet2 B\}=r_{12} (\sheet1 C\otimes\sheet2 B) - (\sheet1 C\otimes\sheet2 B) r_{12}. \label{CB-bracket1}
\eea
It is a standard result that $GL_N\times GL_N$  with the above bracket is a Poisson Lie group.

We classify all quadratic Poisson brackets on $GL_N\times GL_N\times\mathcal A$ for which the Poisson Lie group action of $GL_N\times GL_N$ on $\mathcal A$  defined by (\ref{triple}) is Poisson:

\begin{lm}\label{lem-ABC-new}
Provided that the brackets between $B$ and $C$ matrices are given by (\ref{BB-bracket1}), (\ref{CC-bracket1}), and (\ref{CB-bracket1}),
that the brackets for $\mathbb A$ are given by (\ref{Poisson-r}), and that the brackets between $B$, $C$, and $\mathbb A$
have the form (\ref{bracket-Q}):
\bea
&{}&
\{\sheet1{B}\ocomma\sheet2{\mathbb A}\}=\sheet1{B}Q_{12}\sheet2{\mathbb A}+\sheet1{B}\sheet2{\mathbb A}Q_{12}^{t_2},
\nonumber\\
&{}&
\{\sheet1{C}\ocomma\sheet2{\mathbb A}\}=\sheet1{C}Q_{12}\sheet2{\mathbb A}+\sheet1{C}\sheet2{\mathbb A}Q_{12}^{t_2},
\nonumber
\eea
in order for these brackets to satisfy the Jacobi identities and for
the mapping ${\mathbb A}\mapsto B{\mathbb A}C^{\text{T}}$
to be an automorphism of the Poisson algebra (\ref{Poisson-r}), we have
exactly three choices of the matrix $Q_{12}$ itemized in Lemma~\ref{lm-generalization}.

In cases (ii) and (iii),
the combination ${\mathfrak A}:=B{\mathbb A}^{-\text{T}}C^{\text{T}}$ satisfies the same algebra (\ref{Poisson-r}) as
$\mathbb A$ and $B{\mathbb A}C^{\text{T}}$.
\end{lm}

\proof  Straightforward calculation using the $r$-matrix form of writing for the
corresponding brackets.\endproof

We again have an analogue of Lemma~\ref{duality-ii-iii}.

\begin{lm}\label{duality-ii-iii-ABC}
The transformation $\mathbb A\mapsto \mathbb A^{-1}$, $B\mapsto C^{-\text{T}}$, $C\mapsto B^{-\text{T}}$
is an antiautomorphism of the Poisson algebra for the triple $(\mathbb A,B,C)$ with the
interchange $Q_{12}^{(ii)}\leftrightarrow Q_{12}^{(iii)}$.
\end{lm}

We now interpret $\mathcal A$ as a Poisson symmetric space:

\begin{theorem}\label{th:p-sym}
The map
\begin{equation}\label{map-theta}
\begin{array}{cccc}
\Theta: & GL_N\times GL_N&  \to & GL_N\times GL_N\\
&(B,C)&\mapsto& (C^{-T},B^{-T})
\end{array}
\end{equation}
is an involutive anti-Poisson automorphism so that $(GL_N\times GL_N,\Theta)$ is a Poisson symmetric Lie group.

Let $H\subset GL_N\times GL_N$ be the fixed point set of $\Theta$. Then the immersion:
\begin{equation}
\begin{array}{cccc}
i: & GL_N\times GL \slash H &  \to &\mathcal A\\
&(B,C)&\mapsto&  BC^{\text{T}}
\end{array}
\end{equation}
is a Poisson isomorphism.
\end{theorem}

\proof
The fact that $\Theta$ is an involutive automorphism is obvious. To prove that it is anti-Poisson we need to prove that $(B',C'):= \Theta (B,C)$ satisfy the following Poisson brackets:
\bea
&{}&\{\sheet1{B'}\ocomma\sheet2{B'\}}=- r_{12} (\sheet1{B'}\otimes\sheet2{B'}) + (\sheet1{B'}\otimes\sheet2{B'}) r_{12}, \nn\\
&{}&\{\sheet1{C'}\ocomma\sheet2{C'}\}=-r_{12} (\sheet1{C'}\otimes\sheet2{C'}) + (\sheet1{C'}\otimes\sheet2{C'}) r_{12}, \nn\\
&{}&\{\sheet1{C'}\ocomma\sheet2{B'}\}=- r_{12} (\sheet1{C'}\otimes\sheet2{B'}) + (\sheet1{C'}\otimes\sheet2{B'}) r_{12}.\nn
\eea
The first two are obvious, let us prove the third:
\bea
\{\sheet1{C'}\ocomma\sheet2{B'}\}&=&\{\sheet1{B^{-T}}\ocomma\sheet2{C^{-T}}\}=\sheet1{B^{-T}}\otim\sheet2{C^{-T}} \{\sheet1{B^{\text{T}}}\ocomma\sheet2{C^{\text{T}}}\}  \sheet1{B^{-T}}\otim\sheet2{C^{-T}}=\nn\\
&=&\sheet1{B^{-T}}\otim\sheet2{C^{-T}}  (-r_{21} (\sheet1 B\otimes\sheet2 C)+ (\sheet1 B\otimes\sheet2 C) r_{21})^{\text{T}} \sheet1{B^{-T}}\otim\sheet2{C^{-T}}=\nn\\
&=& -r_{21}^{\text{T}}  \sheet1{B^{-T}}\otim \sheet2{C^{-T}}+\sheet1{B^{-T}}\otim\sheet2{C^{-T}} r_{21}^{\text{T}}=- r_{12} (\sheet1{C'}\otimes\sheet2{B'}) + (\sheet1{C'}\otimes\sheet2{B'}) r_{12},
\nn
\eea
where in the last step we have used that $r_{21}=r_{12}^{\text{T}}$.

Let us now consider the fixed point set $H\subset GL_N\times GL_N$ of $\Theta$:
$$
H=\{(H_1,H_2)\in GL_N\times GL_N| H_1=H_2^{-T}\}.
$$
Then $GL_N\times GL_N\slash H$ is the set of equivalence classes $(B_1,C_1)\sim (B_2,C_2)$ iff $(B_1,C_1)=(B_2 H_1,C_2 H_2)$ for some $ (H_1,H_2)\in H$ and it is straightforward to prove that the immersion $i$ is an isomorphism. To prove that this is a Poisson isomorphism we first observe that the following graph commutes:
$$
\begin{array}{ccc}
GL_N\times GL_N\slash H & \mapi &  \mathcal A\\
\pi\,\big\uparrow \quad \quad& {\nearrow}_{\hat i } &\\
GL_N\times GL_N  & & \\
\end{array}
$$
where $\pi$ is the coset map which associates to each element $(B,C)$ its right coset and $\hat i(B,C)=B C^{\text{T}}$. In  \cite{LuW} it is proved that there exists a unique Poisson bracket on $GL_N\times GL_N\slash H $ for which $\pi$ is a Poisson map. So assuming we endow $GL_N\times GL_N\slash H $ with such a Poisson map, if we prove that $\hat i$ is Poisson than $i$ is. To prove that $\hat i$ is Poisson we just need to prove that
\bea
\{\sheet1{BC^{\text{T}}}\ocomma\sheet2{BC^{\text{T}}}\}&=&r_{12}(\sheet1{BC^{\text{T}}}\otimes\sheet2{BC^{\text{T}}})
-(\sheet1{BC^{\text{T}}}\otimes\sheet2{BC^{\text{T}}})r_{12}+\nn\\
&+&\sheet1{BC^{\text{T}}}r_{12}^{t_1}\sheet2{BC^{\text{T}}}
-\sheet2{BC^{\text{T}}}r_{12}^{t_1}\sheet1{BC^{\text{T}}}.\nn
\eea
This is again a straightforward computation which uses the $r$-matrix properties.
\endproof

\section{Chains of $B$-matrices}\label{s:groupoid-new}
\setcounter{equation}{0}

We now introduce the Poisson structure on the extended space $\mathcal A\otimes_{k=1}^{n} GL_N$  of chains $({\mathbb A},B_1,B_2,\dots,B_n)$. We want to postulate  the brackets between $B_i$ and $B_{i+1}$ that are compatible with the following natural (partial) multiplication operation:
\be
({\mathbb A},B_1)\circ (B_1{\mathbb A}B_1^{\text{T}},B_2)=({\mathbb A},B_2B_1),
\label{group-ABB}
\ee
and its chain analogue
$$
({\mathbb A},B_1)\circ (B_1{\mathbb A}B_1^{\text{T}},B_2)\circ\cdots\circ(B_{j-1}
\cdots B_1{\mathbb A}B_1^{\text{T}}\cdots B_{j-1}^{\text{T}},B_j)=({\mathbb A},B_j B_{j-1}\cdots B_2 B_1).
$$

It is easy to see that if we impose
\bea
&{}&
\{\sheet1{B_2}\ocomma\sheet2{B_1}\}=\sheet1{B_2}Q_{12}\sheet2{B_1},
\nonumber\\
&{}&
\{\sheet1{B_2}\ocomma\sheet2{B_2}\}=r_{12}\sheet1{B_2}\sheet2{B_2}-\sheet1{B_2}\sheet2{B_2}r_{12},
\nonumber\\
&{}&
\{\sheet1{B_2}\ocomma\sheet2{\mathbb A}\}=0,
\nonumber
\eea
where $Q_{12}$ is chosen as in Theorem \ref{lm-generalization},  then these brackets ensure that
all the three mappings $(\mathbb A,B_1,B_2)\mapsto (\mathbb A, B_1)$,
$(\mathbb A,B_1,B_2)\mapsto (\mathbb A, B_2B_1)$, and $(\mathbb A,B_1,B_2)\mapsto (B_1\mathbb AB_1^{\text{T}}, B_2)$
be Poisson.

In the multiple chain generalisation we can prove the following:

\begin{lm}\label{lem-groupoid-Poisson}
The Poisson structure compatible with the groupoid multiple product
$$
({\mathbb A},B_1)\circ (B_1{\mathbb A}B_1^{\text{T}},B_2)\circ\cdots\circ(B_{j-1}
\cdots B_1{\mathbb A}B_1^{\text{T}}\cdots B_{j-1}^{\text{T}},B_j)=({\mathbb A},B_j B_{j-1}\cdots B_2 B_1)
$$
has the following form: the brackets between $\mathbb A$ are given by (\ref{Poisson-r}), the bracket
between $\mathbb A$ and $B_1$ has the form (\ref{bracket-Q}), $B_k$ with $k>2$ Poisson commute with
$\mathbb A$, whereas the Poisson brackets between $B_k$ are
\bea
&{}&
\{\sheet1{B_k}\ocomma\sheet2{B_k}\}=r_{12}\sheet1{B_k}\sheet2{B_k}-\sheet1{B_k}\sheet2{B_k}r_{12},\quad k=1,\dots,j,
\label{BiBi}\\
&{}&
\{\sheet1{B_{k+1}}\ocomma\sheet2{B_k}\}=\sheet1{B_{k+1}}Q^{[k]}_{12}\sheet2{B_k},\quad k=1,\dots,j-1
\label{Bi+1Bi}\\
&{}&
\{\sheet1{B_k}\ocomma\sheet2{B_l}\}=0,\quad |k-l|>1,\nonumber
\eea
where $Q^{[k]}_{12}$ is either zero, or $-r_{12}^{t_2}$, or $r_{12}^{t_1}$. As before, we refer these
three cases to as (i), (ii), and (iii).
\end{lm}

\begin{remark}\label{rem-Qk}
In the relation (\ref{Bi+1Bi}) the $Q^{[k]}$-matrices can be different for different $k$ (of course, all of them must
be of one of three types, (i), (ii), or (iii)). In all calculations below we always assume that $Q^{[k]}=Q^{(ii)}$
for all $k$. In this case, the Poisson relations are uniform, but the structure of central elements is different
for odd and even $j$,  as we shall see in Section \ref{s:central}. If, on the contrary, we set, say, $Q^{[2r+1]}=Q^{(iii)}$ and $Q^{[2r]}=Q^{(ii)}$, we obtain
uniform expressions for central elements for the price of introducing an alternating Poisson brackets.
\end{remark}

\section{Chains of $(B,C)$-pairs}\label{s:groupoid-new-BC}
\setcounter{equation}{0}

We now introduce the
Poisson structure on the extended space $\mathcal A\otimes_{k=1}^{2n} GL_N$
of chains $({\mathbb A},(B_1,C_1),(B_2,C_2),\dots,(B_n,C_n))$.
The (partial) multiplication operation now reads
\be
({\mathbb A},B_1,C_1)\circ (B_1{\mathbb A}C_1^{\text{T}},B_2,C_2)=({\mathbb A},B_2B_1,C_2C_1).
\label{group-ABC}
\ee
We now must postulate the brackets between $B_1,C_1$ and $B_2,C_2$ that are compatible with this multiplication
(\ref{group-ABC}). Inside every pair $(B_i,C_i)$, the brackets coincide with (\ref{BB-bracket1})--(\ref{CB-bracket1}); it is then easy to see that if we impose
\bea
&{}&
\{\sheet1{B_2}\ocomma\sheet2{B_1}\}=\sheet1{B_2}Q_{12}\sheet2{B_1},
\nonumber\\
&{}&
\{\sheet1{C_2}\ocomma\sheet2{B_1}\}=\sheet1{C_2}Q_{12}\sheet2{B_1},
\nonumber\\
&{}&
\{\sheet1{B_2}\ocomma\sheet2{C_1}\}=\sheet1{B_2}Q_{12}\sheet2{C_1},
\nonumber\\
&{}&
\{\sheet1{C_2}\ocomma\sheet2{C_1}\}=\sheet1{C_2}Q_{12}\sheet2{C_1},
\nonumber\\
&{}&
\{\sheet1{B_2}\ocomma\sheet2{\mathbb A}\}=\{\sheet1{C_2}\ocomma\sheet2{\mathbb A}\}=0,
\nonumber
\eea
where $Q_{12}$, as above, is one of three cases (i), (ii), and (iii), then
all the three mappings $(\mathbb A,(B_1,C_1),(B_2,C_2))\mapsto (\mathbb A, B_1,C_1)$,
$(\mathbb A,(B_1,C_1),(B_2,C_2))\mapsto (\mathbb A, B_2B_1,C_2C_1)$, and
$(\mathbb A,(B_1,C_1),(B_2,C_2))\mapsto (B_1\mathbb AC_1^{\text{T}}, B_2,C_2)$
are Poisson.

\begin{lm}\label{lem-groupoid-Poisson-BC}
The Poisson structure compatible with the groupoid multiple product
\bea
&{}&
({\mathbb A},B_1,C_1)\circ (B_1{\mathbb A}C_1^{\text{T}},B_2,C_2)\circ\cdots\circ(B_{j-1}
\cdots B_1{\mathbb A}C_1^{\text{T}}\cdots C_{j-1}^{\text{T}},B_j,C_j)\nonumber\\
&{}&
=({\mathbb A},B_j B_{j-1}\cdots B_2 B_1,C_j C_{j-1}\cdots C_2 C_1)\nonumber
\eea
has the following form: the brackets between $\mathbb A$ are given by (\ref{Poisson-r}), the bracket
between $\mathbb A$ and $B_1$, $C_1$ has the form as in Lemma \ref{lem-ABC-new}, $B_k$ and $C_k$ with $k>2$ Poisson commute with $\mathbb A$, the Poisson brackets between $B_k$ and $C_k$ are of the form
(\ref{BB-bracket1})--(\ref{CB-bracket1}) for any $k$, and the remaining
possibly nonvanishing Poisson brackets are
\bea
&{}&
\{\sheet1{B_{k+1}}\ocomma\sheet2{B_k}\}=\sheet1{B_{k+1}}Q^{[k]}_{12}\sheet2{B_k},\quad k=1,\dots,j-1,
\label{Bk+1Bk}\\
&{}&
\{\sheet1{C_{k+1}}\ocomma\sheet2{B_k}\}=\sheet1{C_{k+1}}Q^{[k]}_{12}\sheet2{B_k},\quad k=1,\dots,j-1,
\label{Ck+1Bk}\\
&{}&
\{\sheet1{B_{k+1}}\ocomma\sheet2{C_k}\}=\sheet1{B_{k+1}}Q^{[k]}_{12}\sheet2{C_k},\quad k=1,\dots,j-1,
\label{Bk+1Ck}\\
&{}&
\{\sheet1{C_{k+1}}\ocomma\sheet2{C_k}\}=\sheet1{C_{k+1}}Q^{[k]}_{12}\sheet2{C_k},\quad k=1,\dots,j-1,
\label{Ck+1Ck}
\eea
where $Q^{[k]}_{12}$ is either zero, or $-r_{12}^{t_2}$, or $r_{12}^{t_1}$. As before, we refer these
three cases to as (i), (ii), and (iii).
\end{lm}

In the above relations, like in the case of
chains of $B$-matrices, the matrices $Q^{[k]}$ can be different for different $k$.
In all calculation below we always assume that $Q^{[k]}=Q^{(ii)}$ for all $k$.

\section{Central elements}\label{s:central}
\setcounter{equation}{0}

\subsection{Central elements of the $\mathbb A$- and $B$-matrix algebras}\label{ss:minors}
The central elements of the Poisson algebra (\ref{Poisson}) of $a_{i,j}$ are of two types: we have
the polynomial central elements and rational central elements; together
they form a set of exactly $N$ algebraically independent central
elements.

The {\it polynomial
central elements} $r_k$ are given
by the coefficients of $\lambda^{-k}$, $k=0,1,\dots, \left[N/2\right]+1$, of the polynomial
$$
\det(\mathbb A +\lambda^{-1}\mathbb A^{\text{T}}).
$$
The {\it rational central elements}  are defined by
the {\it bottom--left minors} of the general matrix $\mathbb A$ (provided these minors do not vanish):
we take
$$
M^-_d(A):=\det\left(\begin{array}{ccc}
a_{N-d+1,1}&\dots&a_{N+d-1,d}\\
\vdots&\dots&\vdots\\
a_{N,1}&\dots&a_{N,d}\end{array}
\right).
$$
In \cite{ChM} (see also \cite{Bon1} where these elements were found to be central for symmetric $\mathbb A$),
we have proved that for every $d=1,\dots,\left[\frac{N}{2}\right]$ the quantities
$$
b_d:=M^-_d(A)/M^-_{N-d}(A)
$$
are central elements of the Poisson algebra (\ref{Poisson}).

The central elements of the Poisson algebra (\ref{Poisson-b}) (see
\cite{FM}) are generated by the complementary minors: let
\be
M^-_d(B):=\det\left(\begin{array}{ccc}
b_{N-d+1,1}&\dots&b_{N+d-1,d}\\
\vdots&\dots&\vdots\\
b_{N,1}&\dots&b_{N,d}\end{array}
\right)
\label{Minor-}
\ee
and
\be
M^+_d(B):=\det\left(\begin{array}{ccc}
b_{1,N-d}&\dots&b_{1,N}\\
\vdots&\dots&\vdots\\
b_{d,N-d}&\dots&b_{d,N}\end{array}
\right)
\label{Minor+}
\ee
be the minors located at the respective
bottom-left and upper-right corners of the matrix $B$. We then have exactly
$N$ algebraically independent central elements
$$
c_d=M^+_d(B)/M^-_{N-d}(B),\quad d=1,\dots, N.
$$
Note that these are exactly the minors that appeared in \cite{Bondal} in the structure of $B$-matrices for the groupoid of
upper-triangular matrices.

\subsection{Casimir functions of the Lie--Poisson brackets for the $B$ and $C$ matrices}\label{ss:Lie-central}

\begin{lm}\label{Lie-Poisson-central}
The Poisson--Lie brackets (\ref{BB-bracket1})--(\ref{CB-bracket1}) of the $(B,C)$-system
possess $2N$ Casimir functions:
\begin{itemize}
\item[(i)] $N+1$ Casimir functions
generated by the minors of the matrices $B$ and $C$ (in the notation of (\ref{Minor-})) and (\ref{Minor+})
\be
M^-_d(B)/M^+_{N-d}(C),\quad d=0,1,\dots,N-1,N;
\ee
\item[(ii)] $N+1$ Casimir functions $q_s$ generated by the coefficients of $\lambda^s$ of the expansion of
\be
\det(B+\lambda C)=\sum_{s=0}^N \lambda^s q_s.
\ee
\end{itemize}
Note that $\det C$ and $\det B$ enter the both sets, so the total number of algebraically independent Casimir
functions is exactly $2N$.
\end{lm}

{\bf Proof.} That the above elements are central is a relatively easy calculation. A more lengthy is the proof that the
general degeneracy of the Poisson brackets (\ref{BB-bracket1})--(\ref{CB-bracket1}) is $2N$. To prove it, we first consider
the linearized algebraic version of these brackets for $B={\mathbb E}+\epsilon b$ and $C={\mathbb E}+\epsilon c$:
\bea
&{}&\{\sheet1 b\ocomma\sheet2 b\}=r_{12} (\sheet1 {\mathbb E}\otimes\sheet2 b+\sheet1 b\otimes\sheet2 {\mathbb E}) -
(\sheet1 {\mathbb E}\otimes\sheet2 b+\sheet1 b\otimes\sheet2 {\mathbb E}) r_{12}, \label{bb-bracket}\\
&{}&\{\sheet1 c\ocomma\sheet2 c\}=r_{12} (\sheet1 {\mathbb E}\otimes\sheet2 c+\sheet1 c\otimes\sheet2 {\mathbb E}) -
(\sheet1 {\mathbb E}\otimes\sheet2 c+\sheet1 c\otimes\sheet2 {\mathbb E}) r_{12}, \label{cc-bracket}\\
&{}&\{\sheet1 c\ocomma\sheet2 b\}=r_{12} (\sheet1 {\mathbb E}\otimes\sheet2 b+\sheet1 c\otimes\sheet2 {\mathbb E}) -
(\sheet1 {\mathbb E}\otimes\sheet2 b+\sheet1 c\otimes\sheet2 {\mathbb E}) r_{12}, \label{cb-bracket}
\eea
We are going to solve the linear system of equations w.r.t. $N\times N$-matrices $x$ and $y$
\be
\{\sheet1 b\ocomma \hbox{tr\,}_2(\sheet2 {bx}+\sheet2 {cy}) \}=\{\sheet1 c\ocomma \hbox{tr\,}_2(\sheet2 {bx}+\sheet2 {cy}) \}=0.
\label{xy}
\ee
Using the explicit form of the $r$-matrix, these two systems of equations can be written in the form (here $P_{+,1/2}$ and $P_{-,1/2}$
are the standard projection operators)
\bea
&{}&P_{+,1/2}(x)b+P_{+,1/2}(bx)-bP_{+,1/2}(x)-P_{+,1/2}(xb)
\nonumber\\
&{}&\qquad -P_{-,1/2}(y)b-P_{-,1/2}(cy)+bP_{-,1/2}(y)+P_{-,1/2}(yc)=0,
\label{eq1}
\\
&{}&P_{+,1/2}(bx)+P_{+,1/2}(x)c-cP_{+,1/2}(x)-P_{+,1/2}(xb)
\nonumber\\
&{}&\qquad +P_{+,1/2}(cy)+P_{+,1/2}(y)c-cP_{+,1/2}(y)-P_{+,1/2}(yc)=0.
\label{eq2}
\eea
Subtracting the second equation from the first one, we obtain a simple restriction that
\be
\bigl[P_{+,1/2}(x)-P_{-,1/2}(y),b-c\bigr]=0.
\label{eq3}
\ee
We now choose the matrices $b$ and $c$ in the special form containing only diagonal and anti-diagonal parts,
$b_{i,j}=b_i \delta_{i,j}+\delta_{i,N+1-i}$ and $c_{i,j}=c_i \delta_{i,j}+\delta_{i,N+1-i}$ with all $b_i$ and $c_i$
distinct. It is then easy to see that among all non-diagonal entries of $x$ and $y$ only the entries on the lower half-anti-diagonal
of $x$ and on the upper half-anti-diagonal of $y$ can be nonzero and substituting this anzatz into (\ref{eq1}) we obtain exactly
$N$ equations
\bea
&{}&x_{d,N+1-d}+y_{N+1-d,d}=0,\quad d=1,\dots [N/2],\nonumber\\
&{}&1/2(x_{d,d}+y_{d,d})=1/2(x_{N+1-d,N+1-d}+y_{N+1-d,N+1-d})\nonumber
\eea
on $3N$ variables. This clearly indicates that we have exactly $2N$-dimensional space of solutions corresponding to
$2N$ Casimir functions. The lemma is proved.\quad$\square$

\subsection{The Casimir functions of the type-(ii) $(\mathbb A,B)$-system}\label{ss:newCasimirs}
We begin constructing Casimir functions for algebra (ii) by noting that the brackets (\ref{Poisson-r}) and
(\ref{bracket-Q}) coincide for ${\mathbb A}$ and ${\mathbb A}^{\text{T}}$, that is,
\be
\{\sheet1{B}\ocomma(\sheet2{\mathbb A}+\lambda \sheet2{\mathbb A}^{\text{T}})\}=
\sheet1{B}Q_{12}(\sheet2{\mathbb A}+\lambda \sheet2{\mathbb A}^{\text{T}})
+\sheet1{B}(\sheet2{\mathbb A}+\lambda \sheet2{\mathbb A}^{\text{T}})Q_{12}^{t_2}
\label{BA-lambda}
\ee
for any $\lambda$ and for any choice of the $r$-matrix $Q$. It then follows, in the case of bracket (ii), that
\be
\{b_{ij},\det({\mathbb A}+\lambda {\mathbb A}^{\text{T}})\}=-2 b_{ij} \det({\mathbb A}+\lambda {\mathbb A}^{\text{T}}),
\label{b-det}
\ee
and, recalling that every $\det({\mathbb A}+\lambda {\mathbb A}^{\text{T}})$ is a Casimir function of the $\mathbb A$-algebra,
we obtain that the elements
\be
\frac{\det({\mathbb A}+\lambda {\mathbb A}^{\text{T}})}{\det {\mathbb A}}
\label{A-central}
\ee
are Casimir functions of the total algebra.

To construct the other set of Casimir functions we recall the combination $\mathfrak A$ introduced
in (\ref{mathfrakA}).

\begin{remark}\label{rem-S}
The matrix ${\mathfrak A}$ has the following Poisson relations with ${\mathbb A}$ and $B$:
\bea
&{}&\{\sheet1{B}\ocomma\sheet2{\mathfrak A}\}=r_{12}\sheet1{B}\sheet2{\mathfrak A}+
\sheet2{\mathfrak A}r_{12}^{t_2}\sheet1{B},\label{B-A-new}\\
&{}&\{\sheet1{\mathfrak A}\ocomma\sheet2{\mathbb A}\}=-2\sheet1{B}
\bigl(P_{12}^{t_1}+\sheet1{\mathbb A}^{-\text{T}}\sheet1{\mathbb A}P_{12}\bigr)\sheet1{B^{\text{T}}}.\label{A-A-new}
\eea
We see that the second bracket destroys the structure of $\mathfrak A$. However, if we introduce
the combination
\be
S:={\mathbb A}^{\text{T}}B^{-1}
\label{S}
\ee
we observe that this quantity does have consistent brackets with ${\mathbb A}$, $B$, and itself:
\bea
&{}&
\{\sheet1{S}\ocomma\sheet2{B}\}=\sheet1{S}r_{12}^{t_1t_2}\sheet2{B}+\sheet2{B}r_{12}^{t_2}\sheet1{S},
\label{S-B}\\
&{}&
\{\sheet1{S}\ocomma\sheet2{\mathbb A}\}=r_{12}\sheet1{S}\sheet2{\mathbb A}+\sheet2{\mathbb A}r_{12}^{t_2}\sheet1{S},
\label{S-A}\\
&{}&
\{\sheet1{S}\ocomma\sheet2{S}\}=r_{12}\sheet1{S}\sheet2{S}-\sheet1{S}\sheet2{S}r_{12}^{t_1t_2}.
\label{S-S}
\eea
These formulas are instrumental when finding Casimir functions for the brackets of type (ii).
\end{remark}

We now use the  algebra (\ref{S-B}--\ref{S-S})
of $S$-variables (\ref{S}) to study the Casimirs. We first
write formulas (\ref{S-B}), (\ref{S-A}),
the Poisson relation between $b_{i,j}$ and $a_{k,l}$, and the brackets between $b$'s and $s$'s in components:
\bea
\{s_{i,j},b_{k,l}\}&=&\sum_{\rho=1}^m s_{i,\rho}b_{\rho,l}\theta(j-\rho)\delta_{j,k}+
s_{l,j}b_{k,i}\theta(i-l),
\label{s-b}\\
\{s_{i,j},a_{k,l}\}&=&s_{k,j}a_{i,l}\theta(i-k)+s_{l,j}a_{k,i}\theta(i-l),
\label{s-a}\\
\{b_{i,j},a_{k,l}\}&=&-b_{i,k}a_{j,l}\theta(k-j)-b_{i,l}a_{k,j}\theta(l-j),
\label{b-a}\\
\{b_{i,j},b_{k,l}\}&=&b_{i,l}b_{k,j}(\theta(i-k)-\theta(l-j)),
\label{b-b}\\
\{s_{i,j},s_{k,l}\}&=&s_{i,l}s_{k,j}(\theta(i-k)-\theta(j-l)),
\label{s-s}
\eea

We now let $M_B^{p}$ denote the $(p\times p)$-minor of the matrix $B$ located at the upper-right corner and
let $M_S^{p}$ denote the $(p\times p)$-principal minor of the matrix $S$ (located at the upper-left corner).

Using the same technique as in Sec.~\ref{ss:minors}, we can demonstrate that
the above relations (\ref{s-b})--(\ref{b-b}) imply that all the brackets of $a_{i,j}$ and $b_{i,j}$ with
$M_S^{p}$ and $M_B^{p}$ are  :
\bea
\{M_S^{p},a_{i,j}\}&=&[D^{p}]_{i,j}\cdot M_S^{p}a_{i,j},\nonumber\\
\{M_S^{p},b_{i,j}\}&=&[D^{p}]_{i,j}\cdot M_S^{p}b_{i,j},\nonumber\\
\{M_B^{p},a_{i,j}\}&=&[F^{p}]_{i,j}\cdot M_B^{p}a_{i,j},\nonumber\\
\{M_B^{p},b_{i,j}\}&=&[G_0^{p}]_{i,j}\cdot M_B^{p}b_{i,j},\nonumber
\eea
where $D^p,F^p,G^p_0$ are integer-valued matrices composed of blocks in which all entries are the same. We represent them graphically as follows:
\be
\label{mat-DFG}
D^p=
\begin{tabular}{r|c|cl|}
\multicolumn{4}{l}{\phantom{XX}p}\\
  \cline{2-4}
  p& 2 & 1 & 1 \\
  \cline{2-4}
&1 & 0 & 0 \\
  &1 & 0 & 0 \\
  \cline{2-4}
\end{tabular}\,,\quad
F^p=
\begin{tabular}{r|rc|c|}
\multicolumn{4}{r}{\phantom{XX}p}\\
  \cline{2-4}
&  0 & 0 & -1 \\
& 0 & 0 & -1 \\
  \cline{2-4}
p&  -1 & -1 & -2 \\
  \cline{2-4}
\end{tabular}\,,\quad
G_0^p=
\begin{tabular}{r|rc|c|}
\multicolumn{4}{r}{\phantom{XX}p}\\
  \cline{2-4}
p&  1 & 1 & 0 \\
  \cline{2-4}
 & 0 & 0 & -1 \\
 & 0 & 0 & -1 \\
  \cline{2-4}
\end{tabular}\,,
\ee
where
we assume that the unit squares are $(p\times p)$-blocks (and that $N-p$ is visually
twice bigger than $p$) and that each block is composed by a matrix with all entries equal to the given integer (for example the $p\times p$ block $2$ is a $p\times p$ matrix with all entries equal to $2$).

From this graphical representation we immediately obtain that the combination $M_S^p/M_B^{N-p}$
has constant brackets with all $a_{i,j}$ for any $p$:
\be
\left\{\frac{M_S^p}{M_B^{N-p}},a_{i,j}\right\}=[D^p-F^{N-p}]_{i,j} \frac{M_S^p}{M_B^{N-p}}a_{i,j}
=2\frac{M_S^p}{M_B^{N-p}}a_{i,j},\quad\forall p,
\ee
This ratio of minors does not still have constant brackets with $b_{i,j}$, but if we take the
{\em product} of two such ratios, then we finally obtain the constant relation for any $i$ and $j$:
$$
\left\{\frac{M_S^pM_S^{N-p}}{M_B^{N-p}M_B^p},b_{i,j}\right\}=[D^p+D^{N-p}-G_0^p-G_0^{N-p}]_{i,j}
\frac{M_S^pM_S^{N-p}}{M_B^{N-p}M_B^p}b_{i,j}=2\frac{M_S^pM_S^{N-p}}{M_B^{N-p}M_B^p}b_{i,j}.
$$
Now in order to obtain commuting quantities it suffices to multiply this expression by the appropriate
powers of the determinants $\det {\mathbb A}$ and $\det B$. We have therefore proved the following theorem.

\begin{theorem}\label{th-central-A-B}
We have the following $2\left[\frac{N}{2}\right]$
(algebraically independent) Casimir functions of the Poisson brackets (ii):
\begin{itemize}
\item $\left[\frac{N}{2}\right]$ coefficients $Y_p$ of $\lambda$-expansion
$$
\frac{\det(\lambda {\mathbb A}+\lambda^{-1}{\mathbb A}^{\text{T}})}{\det {\mathbb A}}
=(\lambda^N+\lambda^{-N})+\sum_{p=1}^{[N/2]}
(\lambda^{N-2p}+\lambda^{2p-N})Y_p;
$$
\item $\left[\frac{N}{2}\right]$ Casimir functions
$$
X_p:=\frac{M_S^pM_S^{N-p}}{M_B^{N-p}M_B^p}\cdot\frac{[\det B]^2}{\det{\mathbb A}}, \quad p=1,\dots,[N/2],
$$
where $M_S^q$ are the $(q\times q)$-principal (situated at the upper-left corner)
minors of the matrix $S={\mathbb A}^{\text{T}} B^{-1}$ and $M_B^q$ are $(q\times q)$-minors of the matrix $B$
situated at the upper-right corner.
\end{itemize}
\end{theorem}

\begin{remark}\label{rem-homoneneity}
Note that both $Y_p$ and $X_p$ are   functions of $a_{i,j}$ and $b_{i,j}$ of degree zero:
scaling independently all $b_{i,j}$ and all $a_{i,j}$ does not change the values of the Casimir functions.
\end{remark}

\subsection{The Casimir functions of the type-(ii) $(\mathbb A, B, C)$-triple}\label{ss:newABC-Casimirs}

It is again useful to consider the same combination $S={\mathbb A}^{\text{T}}B^{-1}$ as in the previous section. We now need
its commutation relations with $C$:
\be
\{\sheet1{S}\ocomma\sheet2{C}\}=\sheet1{S}r_{12}^{t_1t_2}\sheet2{C}+\sheet2{C}r_{12}^{t_2}\sheet1{S},
\label{S-C}
\ee
We see that they have exactly the same form as (\ref{S-B}), so it is not surprising that the principal minors of the
$S$-matrix have   commutation relations with entries of the $C$ matrix as well. On the other hand, it is
now minors of the matrix $C$, not $B$, located at the upper-right corner that have   commutation relations. 
Applying the same technique as in the proof of Theorem~\ref{th-central-A-B}, we come to the following statement.

\begin{theorem}\label{th-central-A-B-C}
We have the following $2\left[\frac{N}{2}\right]+N$
(algebraically independent) Casimir functions of the Poisson brackets (ii) in the case of $(\mathbb A,B,C)$-triple:
\begin{itemize}
\item $\left[\frac{N}{2}\right]$ coefficients $Y_p$ of $\lambda$-expansion
$$
\frac{\det(\lambda {\mathbb A}+\lambda^{-1}{\mathbb A}^{\text{T}})}{\det{\mathbb A}}
=(\lambda^N+\lambda^{-N})+\sum_{p=1}^{[N/2]}
(\lambda^{N-2p}+\lambda^{2p-N})Y_p;
$$
\item $\left[\frac{N}{2}\right]$ Casimir functions
$$
X_p:=\frac{M_S^pM_S^{N-p}}{M_C^{N-p}M_C^p}\cdot\frac{\det B\det C}{\det{\mathbb A}}, \quad p=1,\dots,[N/2],
$$
where $M_S^q$ are the $(q\times q)$-principal (situated at the upper-left corner)
minors of the matrix $S={\mathbb A}^{\text{T}} B^{-1}$ and $M_C^q$ are $(q\times q)$-minors of the matrix $C$
situated at the upper-right corner;
\item $N$ coefficients $Z_p$ of the $\lambda$-expansion
$$
\det (B+\lambda C)/\det B=1+\sum_{p=1}^N \lambda^p Z_p.
$$
\end{itemize}
\end{theorem}

\subsection{Casimir functions for chains of $B$-matrices}\label{ss:Casimirs4chains}

As above, we mainly consider the case (ii). In this case, for a chain of matrices $B_1,\dots,B_j$ we again have
the special matrix $S$ whose form will be different for odd and even $j$:
\be
\label{S-string}
S=\left\{ \begin{array}{ll}
            {\mathbb A}^{\text{T}} B_1^{-1}B_2^{\text{T}}B_3^{-1}\cdots B_{j-1}^{\text{T}}B_j^{-1} & \hbox{for odd $j$}, \\
            {\mathbb A}^{\text{T}} B_1^{-1}B_2^{\text{T}}B_3^{-1}\cdots B_{j-1}^{-1}B_j^{\text{T}}  & \hbox{for even $j$}.
          \end{array}
\right.
\ee
We now consider the case $j>1$.
The matrix $S$ has the same Poisson relations with ${\mathbb A}$ and $B_1$ irrespectively whether
$j$ is odd or even:
\bea
&{}&
\{\sheet1{S}\ocomma\sheet2{\mathbb A}\}=r_{12}\sheet1{S}\sheet2{\mathbb A}
+\sheet2{\mathbb A}r_{12}^{t_2}\sheet1{S}
\label{AS-string}\\
&{}&
\{\sheet1{S}\ocomma\sheet2{B_1}\}=\sheet2{B_1}r_{12}^{t_2}\sheet1{S}
\label{B1S-string}
\eea
Next, it is easy to see that $S$ commutes with all $B_k$, $k=2,\dots,j-1$. This follows from the following
observation:
\bea
&{}&
\{\sheet1{B_k}\ocomma\sheet2{{B^{\text{T}}}_{k-1}B^{-1}_k{B^{\text{T}}}_{k+1}}\}= \{\sheet1{B_k}\ocomma\sheet2{{B^{\text{T}}}_{k-1}}\}\sheet2{B^{-1}_{k}}\sheet2{{B^{\text{T}}}_{k+1}}
\nonumber\\
&{}&\quad+\sheet2{{B^{\text{T}}}_{k-1}} \{\sheet1{B_k}\ocomma\sheet2{B^{-1}_{k}}\}\sheet2{{B^{\text{T}}}_{k+1}}+
\sheet2{{B^{\text{T}}}_{k-1}} \sheet2{B^{-1}_{k}} \{\sheet1{B_k}\ocomma\sheet2{{B^{\text{T}}}_{k+1}}\}
\nonumber\\
&{}&\quad
=-\sheet1{B_k}\sheet2{{B^{\text{T}}}_{k-1}}r_{12}\sheet2{B^{-1}_{k}}\sheet2{{B^{\text{T}}}_{k+1}}+
\sheet2{{B^{\text{T}}}_{k-1}}\left(-\sheet2{B^{-1}_{k}}r_{12}\sheet1{B_k}+\sheet1{B_k} r_{12}\sheet2{B^{-1}_{k}} \right)\sheet2{{B^{\text{T}}}_{k+1}}\nonumber\\
&{}&\qquad\quad+ \sheet2{{B^{\text{T}}}_{k-1}}\sheet2{B^{-1}_{k}}r_{12} \sheet1{B_k}\sheet2{{B^{\text{T}}}_{k+1}}
\nonumber\\
&{}&\quad =0.\nonumber
\eea
The only difference occurs in the last commutation relations:
\begin{itemize}
\item {\bf for $j$ odd} we have that
$$
\{\sheet1{S}\ocomma\sheet2{B^{-1}_j}\}=-\sheet1{S}\sheet2{B^{-1}_j}r_{12}^{t_1t_2}
$$
and
$$
\{\sheet1{S}\ocomma\sheet2{S}\}=r_{12}\sheet1{S}\sheet2{S}-\sheet1{S}\sheet2{S}r_{12}^{t_1t_2}
$$
\item {\bf for $j$ even} we have that
$$
\{\sheet1{S}\ocomma\sheet2{{B^{\text{T}}}_j}\}=-\sheet1{S}\sheet2{{B^{\text{T}}}_j}r_{12}
$$
and
$$
\{\sheet1{S}\ocomma\sheet2{S}\}=r_{12}\sheet1{S}\sheet2{S}-\sheet1{S}\sheet2{S}r_{12}
$$
\end{itemize}

As in Sec.~\ref{ss:newCasimirs}, we can construct commuting elements from the corresponding minors of $S$ and
$B_k$. But now we again observe the difference between the cases of odd and even $j$.

\begin{defin}
We let $M_S^p$ denote the {\em principal} (located at the upper-left corner) $(p\times p)$-minor
of the matrix $S$ (\ref{S-string}) for {\em odd} $j$ and let the same symbol
$M_S^p$ denote the {\em upper-right} $(p\times p)$-minor of the matrix $S$ (\ref{S-string})
for {\em even} $j$. For all the matrices $B_k$ we let $M_{B_k}^p$ denote their
{\em upper-right} $(p\times p)$-minors.
\end{defin}

As above, all these minors have   commutation relations with all $a_{s,q}$ and $[b_k]_{s,q}$. Below we
present all nonzero commutation relations:
\bea
&{}&\{M_S^p, a_{s,q}\}=[D^p]_{s,q}M_S^p a_{s,q},\nonumber\\
&{}&\{M_S^p, [b_1]_{s,q}\}=[G_-^p]_{s,q}M_S^p [b_1]_{s,q},\nonumber\\
&{}&\{M_S^p, [b_j]_{s,q}\}=\left\{\begin{array}{ll}
                                    [G_+^p]_{s,q}M_S^p [b_j]_{s,q} & \hbox{for even $j$}, \cr
                                    [E^p]_{s,q}M_S^p [b_j]_{s,q} & \hbox{for odd $j$},
                                  \end{array}
\right.\nonumber\\
&{}&\{M_{B_1}^p, a_{s,q}\}=[F^p]_{s,q}M_{B_1}^p a_{s,q},\nonumber\\
&{}&\{M_{B_k}^p, [b_k]_{s,q}\}=[G_0^p]_{s,q}M_{B_k}^p [b_k]_{s,q}, \quad k=1,\dots,j,\nonumber\\
&{}&\{M_{B_k}^p, [b_{k+1}]_{s,q}\}=[G_-^p]_{s,q}M_{B_k}^p [b_{k+1}]_{s,q}, \quad k=1,\dots,j-1,\nonumber\\
&{}&\{M_{B_k}^p, [b_{k-1}]_{s,q}\}=[G_+^p]_{s,q}M_{B_k}^p [b_{k-1}]_{s,q}, \quad k=2,\dots,j.\nonumber
\eea
Here the matrices $D^p$, $F^p$, and $G_0^p$ are defined in (\ref{mat-DFG}), and the remaining three matrices
have the form
\be
\label{mat-EGG}
G_-^p=
\begin{tabular}{|c|cl|}
\multicolumn{3}{l}{p}\\
  \hline
  1 & 0 & 0 \\
  1 & 0 & 0 \\
  1 & 0 & 0 \\
  \hline
\end{tabular}\,,\quad
G_+^p=
\begin{tabular}{r|rcc|}
  \cline{2-4}
&  0 & 0 & 0 \\
 & 0 & 0 & 0 \\
  \cline{2-4}
p&  -1 & -1 & -1 \\
  \cline{2-4}
\end{tabular}\,,\quad
E^p=
\begin{tabular}{r|rcc|}
  \cline{2-4}
p&  1 & 1 & 1 \\
  \cline{2-4}
&  0 & 0 & 0 \\
 & 0 & 0 & 0 \\
  \cline{2-4}
\end{tabular}\,.
\ee

Introducing the combination
\be
\label{Kp}
K^p=\left\{\begin{array}{ll}
\frac{M_S^p}{M_{B_1}^{N-p}}\frac{M_{B_2}^p}{M_{B_3}^{N-p}}\cdots \frac{M_{B_{j-2}}^p}{M_{B_{j-1}}^{N-p}}M_{B_j}^p
& \hbox{for even $j$}, \\
\frac{M_S^p}{M_{B_1}^{N-p}}\frac{M_{B_2}^p}{M_{B_3}^{N-p}}\cdots \frac{M_{B_{j-1}}^p}{M_{B_{j}}^{N-p}}
                                    & \hbox{for odd $j$},
                                  \end{array}
\right.
\ee
we see that it has zero brackets with all $B_k$, $k=1,\dots,j-1$, it has the constant bracket
$$
\{K^p,a_{s,p}\}=2K^p a_{s,p}
$$
with entries of the matrix ${\mathbb A}$, whereas its brackets with $B_j$ are again different depending on
case of which $j$, odd or even, we are dealing with:
\be
\label{Kp-Bj}
\{K^p,[b_j]_{s,p}\}=\left\{\begin{array}{ll}
[G^p_+-G_-^{N-p}+G_0^p]_{s,p}K^p[b_j]_{s,p}
& \hbox{for even $j$}, \cr
[E^{p} +G_-^{p}-G_0^{N-p}]_{s,p}K^p[b_j]_{s,p}
                                    & \hbox{for odd $j$}.
                                  \end{array}
\right.
\ee

In the case of even $j$, the matrix combination $G^p_+-G_-^{N-p}+G_0^p$ has the form
$$
G^p_+-G_-^{N-p}+G_0^p=
\begin{tabular}{r|rcc|}
  \cline{2-4}
p&  0 & 0 & 0 \\
  \cline{2-4}
&  -1 & -1 & -1 \\
  \cline{2-4}
p&  -2 & -2 & -2 \\
  \cline{2-4}
\end{tabular}\,
$$
for $p<\bigl[\frac m2\bigr]$, and it has exactly the same form if we replace $p$ by $N-p$, so
$$
[G^p_+-G_-^{N-p}+G_0^p]-[G^{N-p}_+-G_-^{p}+G_0^{m-p}]=0,
$$
and we obtain that the ratio $K^p/K^{N-p}$ is truly central for any $p=0,\dots,N$. Note that the total number of
such algebraically independent combinations is $[(N+1)/2]$, so, together with $[N/2]$ Casimir functions
generated by $\det({\mathbb A}+\lambda {\mathbb A}^{\text{T}})/\det {\mathbb A}$, we obtain exactly $N$ algebraically
independent Casimir functions.

In the case of odd $j$, the combination $E^{p} +G_-^{p}-G_0^{N-p}$ has the form
$
\begin{tabular}{r|rcc|}
  \cline{2-4}
p&  1 & 1 & 1 \\
  \cline{2-4}
&  0 & 0 & 0 \\
  \cline{2-4}
p&  1 & 1 & 1 \\
  \cline{2-4}
\end{tabular}\,
$ (the upper and lower strips have widths equal $p$) in the case where $p<[N/2]$ and the form
$
\begin{tabular}{r|rcc|}
  \cline{2-4}
N{-}p&  1 & 1 & 1 \\
  \cline{2-4}
&  2 & 2 & 2 \\
  \cline{2-4}
N{-}p&  1 & 1 & 1 \\
  \cline{2-4}
\end{tabular}\,
$ (the upper and lower strips have widths equal $N-p$) in the case where $p>[N/2]$. It is now the {\em product}
$K^pK^{N-p}$, not the ratio, that has constant brackets with all $a$'s and $b_k$'s:
$$
\{K^pK^{N-p},a_{s,p}\}=4 K^pK^{N-p}a_{s,p},\qquad \{K^pK^{N-p},[b_j]_{s,p}\}=2 K^pK^{N-p} [b_j]_{s,p},
$$
and all other brackets vanish. We again (as in Sec.~\ref{ss:newCasimirs}) must multiply by the proper
combination of determinants of ${\mathbb A}$ and $B_k$ (all these determinants have constant brackets with
all $a$'s and $b_k$'s) in order to obtain the true Casimir functions. Their number is $[N/2]$.
We collect together the obtained results in the following theorem.

\begin{theorem}\label{th-central-A-B-chain}
We have the following
(algebraically independent) Casimir functions of the Poisson brackets (ii) in the case of
the chain of matrices ${\mathbb A},\,B_1,\, B_2, \dots, B_j$:
\begin{itemize}
\item {\bf for all $j$} we have $\left[\frac{N}{2}\right]$ coefficients $Y_p$ of $\lambda$-expansion
$$
\frac{\det(\lambda {\mathbb A}+\lambda^{-1}{\mathbb A}^{\text{T}})}{\det{\mathbb A}}
=(\lambda^N+\lambda^{-N})+\sum_{p=1}^{[N/2]}(\lambda^{N-2p}+\lambda^{2p-N})Y_p;
$$
\item {\bf in the case of even $j$} we have $\left[\frac{N+1}{2}\right]$ Casimir functions
$$
X_p=\frac{M^p_S\prod_{k=1}^j M_{B_k}^p}{M^{N-p}_S\prod_{k=1}^j M_{B_k}^{N-p}}
$$
where $M_{B_k}^q$ are $(q\times q)$-upper-right minors of the matrices $B_k$ and $M_S^q$ are the $(q\times q)$-upper-right
minors of the matrix $S={\mathbb A}^{\text{T}} B_1^{-1}B_2^{\text{T}}\cdots B_{j-1}^{-1}B_j^{\text{T}}$;
\item {\bf in the case of odd $j=2r+1$} we have $\left[\frac{N}{2}\right]$ Casimir functions
$$
X_p:=\frac{M_S^pM_S^{N-p}\prod_{k=1}^r\left(M_{B_{2k}}^p M_{B_{2k}}^{N-p}\right)}
{\prod_{k=1}^{r+1}\left(M_{B_{2k-1}}^p M_{B_{2k-1}}^{N-p}\right)}\cdot
\frac{\prod_{k=1}^{r+1}[\det B_{2k-1}]^2}{\det {\mathbb A} \prod_{k=1}^{r}[\det B_{2k}]^2},
$$
where $M_{B_k}^q$ are $(q\times q)$-upper-right minors of the matrices $B_k$ and $M_S^q$ are now the $(q\times q)$-principal (situated at the upper-left corner)
minors of the matrix $S={\mathbb A}^{\text{T}} B_1^{-1}B_2^{\text{T}}\cdots B_{j-1}^{\text{T}}B_j^{-1}$.
\end{itemize}
We therefore have $[\frac N2]+[\frac{N+1}2]=N$ Casimir functions for even $j$ and
$2[\frac N2]$ Casimir functions for odd $j$.
\end{theorem}

\begin{remark}
Note that the statement of Theorem \ref{th-central-A-B-chain} remains valid both for $j=0$ and $j=1$.
In the first case, we obtain all the Casimir functions of the ${\mathbb A}$-matrix algebra (\ref{Poisson-r}) and
in the second case we reproduce the statement of Theorem~\ref{th-central-A-B}.
\end{remark}

\subsection{Casimir functions for chains of $(B,C)$-pairs}\label{ss:Casimirs4BCchains}

As above, we consider the case (ii). In this case, for a chain of $(B,C)$-pairs $(B_1,C_1),\dots,(B_j,C_j)$
we again have the special matrix $S$ whose form is different for odd and even $j$; this matrix has
exactly the form (\ref{S-string}), that is, it is {\em independent} on matrices $C_k$.

The subsequent reasonings are similar to those presented in preceding subsections; we formulate only the final
statement.

\begin{defin}
We let $M_S^p$ denote the {\em principal} (located at the upper-left corner) $(p\times p)$-minor
of the matrix $S$ (\ref{S-string}) for {\em odd} $j$ and let the same symbol
$M_S^p$ denote the {\em upper-right} $(p\times p)$-minor of the matrix $S$ (\ref{S-string})
for {\em even} $j$. For all the matrices $B_k$ we let $M_{B_k}^p$ denote their
{\em upper-right} $(p\times p)$-minors.
\end{defin}

\begin{theorem}\label{th-central-A-BC-chain}
The following are
(algebraically independent) Casimir functions of the Poisson brackets (ii) in the case of
the chain of matrices ${\mathbb A},\,(B_1,C_1),\, \dots, (B_j,C_j)$:
\begin{itemize}
\item {\bf for all $j$} we have $\left[\frac{N}{2}\right]$ coefficients $Y_p$ of $\lambda$-expansion
$$
\frac{\det(\lambda {\mathbb A}+\lambda^{-1}{\mathbb A}^{\text{T}})}{\det{\mathbb A}}
=(\lambda^N+\lambda^{-N})+\sum_{p=1}^{[N/2]}(\lambda^{N-2p}+\lambda^{2p-N})Y_p;
$$
\item {\bf for all $j$} we have $N\cdot j$ coefficients $Z^{(k)}_p$ of $\lambda$-expansion
$$
\frac{\det(B_k+\lambda C_k)}{\det B_k}
=1+\sum_{p=1}^{N}\lambda^{p} Z_p^{(k)}, \quad k=1,\dots,j;
$$
\item {\bf in the case of even $j$} we have $\left[\frac{N+1}{2}\right]$ Casimir functions
$$
X_p=\frac{M^p_S\prod_{k=1}^j M_{C_k}^p}{M^{N-p}_S\prod_{k=1}^j M_{C_k}^{N-p}}
$$
where $M_{C_k}^q$ are $(q\times q)$-upper-right minors of the matrices $C_k$ and $M_S^q$ are the $(q\times q)$-upper-right
minors of the matrix $S={\mathbb A}^{\text{T}} B_1^{-1}B_2^{\text{T}}\cdots B_{j-1}^{-1}B_j^{\text{T}}$;
\item {\bf in the case of odd $j=2r+1$} we have $\left[\frac{N}{2}\right]$ Casimir functions
$$
X_p:=\frac{M_S^pM_S^{N-p}\prod_{k=1}^r\left(M_{C_{2k}}^p M_{C_{2k}}^{N-p}\right)}
{\prod_{k=1}^{r+1}\left(M_{C_{2k-1}}^p M_{C_{2k-1}}^{N-p}\right)}\cdot
\frac{\prod_{k=1}^{r+1}[\det B_{2k-1}\det C_{2k-1}]}{\det {\mathbb A} \prod_{k=1}^{r}[\det B_{2k}\det C_{2k}]},
$$
where $M_S^q$ are now the $(q\times q)$-principal (situated at the upper-left corner)
minors of the matrix $S={\mathbb A}^{\text{T}} B_1^{-1}B_2^{\text{T}}\cdots B_{j-1}^{\text{T}}B_j^{-1}$
and $M_{C_k}^q$ are $(q\times q)$-upper-right minors of the matrices $C_k$.
\end{itemize}
We therefore have $[\frac N2]+[\frac{N+1}2]+Nj=N(j+1)$ Casimir functions for even $j$ and
$2[\frac N2]+Nj$ Casimir functions for odd $j$; because the total number of matrix entries
is $N^2(2j+1)$ it is easy to see that the highest dimension of Poisson leaves is always even.
\end{theorem}

As in the case of the $B$-matrix chains, the statement of Theorem \ref{th-central-A-BC-chain} remains
valid both for $j=0$ and $j=1$.

\begin{remark}
In the case (iii), the statement of Proposition~\ref{prop-anti} remains valid as stated, the analogue of the matrix
$S$ is $\mathfrak S:={\mathbb A}B^{-1}$, and this matrix has the following Poisson relations:
\bea
&{}&
\{\sheet1{\mathfrak S}\ocomma\sheet2{B}\}=
-\sheet1{\mathfrak S}r_{12}\sheet2{B}-\sheet2{B}r_{12}^{t_1}\sheet1{\mathfrak S},
\label{S-B-iii}\\
&{}&
\{\sheet1{\mathfrak S}\ocomma\sheet2{\mathbb A}\}
=-r_{12}^{t_1t_2}\sheet1{\mathfrak S}\sheet2{\mathbb A}-\sheet2{\mathbb A}r_{12}^{t_1}\sheet1{\mathfrak S},
\label{S-A-iii}\\
&{}&
\{\sheet1{\mathfrak S}\ocomma\sheet2{\mathfrak S}\}
=r_{12}\sheet1{\mathfrak S}\sheet2{\mathfrak S}-\sheet1{\mathfrak S}\sheet2{\mathfrak S}r_{12}^{t_1t_2}.
\label{S-S-iii}
\eea
The statement of Theorem~\ref{th-central-A-B-chain} also remains valid provided we make the following changes:
we replace the matrix $S$ by the matrix $\mathfrak S$, which has the form
$\mathfrak S={\mathbb A}B_1^{-1}B_2^{\text{T}}\cdots B_{j-1}^{-1}B_j^{\text{T}}$ for even $j$ and
$\mathfrak S={\mathbb A}B_1^{-1}B_2^{\text{T}}\cdots B_{j-1}^{\text{T}}B_j^{-1}$ for odd $j$, the minors
$M_{\mathfrak S}^q$ are now {\em bottom-right} principal minors for odd $j$ (including $j=1$)
and are {\em lower-left} minors
for even $j$ (including $j=0$), whereas all the minors $M_{B_k}^q$ of all the matrices $B_k$ are {\em lower-left} minors.

For the {\em triple} ${\mathbb A}, B, C$, the central
elements $X_p$ in Theorem~\ref{th-central-A-B-C} have now the form
$$
X_p:=\frac{M_{\mathfrak S}^pM_{\mathfrak S}^{N-p}}{M_B^{N-p}M_B^p}
\cdot\frac{\det B\det C}{\det{\mathbb A}}, \quad p=1,\dots,[N/2],
$$
where $M_{\mathfrak S}^q$ are the {\em bottom-right} minors of the new matrix ${\mathfrak S}:={\mathbb A}C^{-1}$,
and $M_B^q$ are the {\em lower-left} minors of the matrix $B$.
\end{remark}

\section{Poisson and Dirac reductions}\label{s:reductions}
\setcounter{equation}{0}

\subsection{Poisson reductions of the $B$-matrix Poisson algebra}

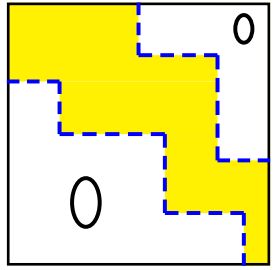
\begin{figure}[tb]
{\psset{unit=0.7}
\begin{pspicture}(-2.5,-2.5)(2.5,2.5)
\psframe[linecolor=yellow, fillstyle=solid, fillcolor=yellow](-2.5,1)(2.5,2.5)
\psframe[linecolor=yellow, fillstyle=solid, fillcolor=yellow](-1.5,0)(2.5,1)
\psframe[linecolor=yellow, fillstyle=solid, fillcolor=yellow](0.5,-1.5)(2.5,0)
\psframe[linecolor=yellow, fillstyle=solid, fillcolor=yellow](2,-2.5)(2.5,-1.5)
\psframe[linecolor=white, fillstyle=solid, fillcolor=white](0,1.5)(2.5,2.5)
\psframe[linecolor=white, fillstyle=solid, fillcolor=white](1.5,-.5)(2.5,2.5)
\psframe[linewidth=1pt](-2.5,-2.5)(2.5,2.5)
\pcline[linecolor=blue, linestyle=dashed, linewidth=1.5pt](-2.5,1)(-1.5,1)
\pcline[linecolor=blue, linestyle=dashed, linewidth=1.5pt](-1.5,1)(-1.5,0)
\pcline[linecolor=blue, linestyle=dashed, linewidth=1.5pt](-1.5,0)(.5,0)
\pcline[linecolor=blue, linestyle=dashed, linewidth=1.5pt](.5,0)(.5,-1.5)
\pcline[linecolor=blue, linestyle=dashed, linewidth=1.5pt](.5,-1.5)(2,-1.5)
\pcline[linecolor=blue, linestyle=dashed, linewidth=1.5pt](2,-1.5)(2,-2.5)
\pcline[linecolor=blue, linestyle=dashed, linewidth=1.5pt](0,2.5)(0,1.5)
\pcline[linecolor=blue, linestyle=dashed, linewidth=1.5pt](0,1.5)(1.5,1.5)
\pcline[linecolor=blue, linestyle=dashed, linewidth=1.5pt](1.5,1.5)(1.5,-.5)
\pcline[linecolor=blue, linestyle=dashed, linewidth=1.5pt](1.5,-.5)(2.5,-.5)
\psellipse[linewidth=1.5pt](-1,-1.3)(0.3,0.5)
\psellipse[linewidth=1.5pt](2,2)(0.2,0.3)
\end{pspicture}
}
\caption{A general Poisson reduction of the algebra (\ref{Poisson-b}). All the items
below the lower dashed broken line and above the upper dashed broken line are zeros.}
\label{fig:Young}
\end{figure}

\begin{lm}\label{lm:B-reduction}
Any reduction depicted in
Fig.~\ref{fig:Young} where all elements below the lower broken line that
goes as in the figure and all elements above the second broken line are set to be zeros
is a Poisson reduction of algebra (\ref{Poisson-b}).
\end{lm}

\subsection{Reductions of the Lie--Poisson brackets for the $B$ and $C$ matrices}\label{ss:Lie-reduction}

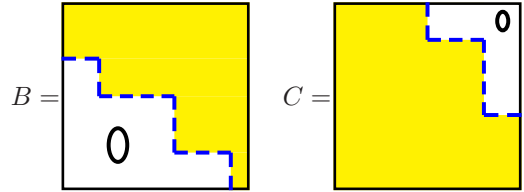
\begin{figure}[tb]
{\psset{unit=0.5}
\begin{pspicture}(-3.5,-2.5)(3.5,2.5)
\rput(-3.2,0){\makebox(0,0){$B=$}}
\psframe[linecolor=yellow, fillstyle=solid, fillcolor=yellow](-2.5,1)(2.5,2.5)
\psframe[linecolor=yellow, fillstyle=solid, fillcolor=yellow](-1.5,0)(2.5,1)
\psframe[linecolor=yellow, fillstyle=solid, fillcolor=yellow](0.5,-1.5)(2.5,0)
\psframe[linecolor=yellow, fillstyle=solid, fillcolor=yellow](2,-2.5)(2.5,-1.5)
\psframe[linewidth=1pt](-2.5,-2.5)(2.5,2.5)
\pcline[linecolor=blue, linestyle=dashed, linewidth=1.5pt](-2.5,1)(-1.5,1)
\pcline[linecolor=blue, linestyle=dashed, linewidth=1.5pt](-1.5,1)(-1.5,0)
\pcline[linecolor=blue, linestyle=dashed, linewidth=1.5pt](-1.5,0)(.5,0)
\pcline[linecolor=blue, linestyle=dashed, linewidth=1.5pt](.5,0)(.5,-1.5)
\pcline[linecolor=blue, linestyle=dashed, linewidth=1.5pt](.5,-1.5)(2,-1.5)
\pcline[linecolor=blue, linestyle=dashed, linewidth=1.5pt](2,-1.5)(2,-2.5)
%
\psellipse[linewidth=1.5pt](-1,-1.3)(0.3,0.5)
\end{pspicture}
\begin{pspicture}(-3.5,-2.5)(3.5,2.5)
\rput(-3.2,0){\makebox(0,0){$C=$}}
\psframe[linecolor=yellow, fillstyle=solid, fillcolor=yellow](-2.5,-2.5)(2.5,2.5)
\psframe[linecolor=white, fillstyle=solid, fillcolor=white](0,1.5)(2.5,2.5)
\psframe[linecolor=white, fillstyle=solid, fillcolor=white](1.5,-.5)(2.5,2.5)
\psframe[linewidth=1pt](-2.5,-2.5)(2.5,2.5)
\pcline[linecolor=blue, linestyle=dashed, linewidth=1.5pt](0,2.5)(0,1.5)
\pcline[linecolor=blue, linestyle=dashed, linewidth=1.5pt](0,1.5)(1.5,1.5)
\pcline[linecolor=blue, linestyle=dashed, linewidth=1.5pt](1.5,1.5)(1.5,-.5)
\pcline[linecolor=blue, linestyle=dashed, linewidth=1.5pt](1.5,-.5)(2.5,-.5)
\psellipse[linewidth=1.5pt](2,2)(0.2,0.3)
\end{pspicture}
}
\caption{A general Poisson reduction of the $B$--$C$ algebra (\ref{BB-bracket1})-(\ref{CB-bracket1}). All the items
below the lower dashed broken line for the matrix $B$ and above the upper dashed broken line for the matrix $C$ are zeros.}
\label{fig:Young2}
\end{figure}

\begin{lm}\label{lm:CB-reduction}
Any reduction depicted in
Fig.~\ref{fig:Young2} where all elements of the matrix $B$ that are below the broken line
and all elements of the matrix $C$ that are above the  broken line are set to be zeros
is a Poisson reduction of algebra (\ref{BB-bracket1})-(\ref{CB-bracket1}).
\end{lm}

\begin{lm}\label{lm:CB-B-reduction}
The reduction of the $(\mathbb A, B, C)-triple$ to the $(\mathbb A,B)$-double realized by
imposing the constraint $B=C$ is Poisson for any choice of the matrix $Q_{12}$.
\end{lm}

\begin{remark}\label{rem-BC-red}
Note that the above Poisson reductions of $B$ and $C$ matrices are compatible with Poisson reductions of
${\mathbb A}$-algebra (\ref{Poisson}): in case of Lemma~\ref{lm:CB-reduction} if the matrices
$B$ and $C^{\text{T}}$ have a b.u.t. form, then the matrix $\mathbb A=BC^{\text{T}}$ will have the same b.u.t. form, which
is known to be
a Poisson reduction of the algebra (\ref{Poisson}); in case of Lemma~\ref{lm:CB-B-reduction}, choosing
$B=C$ corresponds to restricting $\mathbb A$ to a symmetric form, which is again a Poisson reduction of
algebra (\ref{Poisson}).
\end{remark}

\subsection{Restriction to the groupoid of b.u.t. matrices}\label{ss:Dirac}
We now implement the standard procedure of Dirac reduction \cite{Dirac} (which is of extensive use
in quantum field theory, see, e.g., \cite{Wein}) for constructing brackets on manifolds  defined by constraints of the second kind.

Let us recall this procedure here: given a Poisson manifold of dimension $M$,  denote the coordinates by  $X_i$, $i=1,\dots,M$. Choose a set of  $2N$ constraints $C_k(X)$ of the second kind (i.e. such that their mutual Poisson brackets restricted to the constraint manifold are not all zero), if the matrix  $D_{k,l}(X):=\{C_k,C_l\}$ is
nondegenerate on the constraint surface $C_k(X)=0$, $k=1,\dots,2N$, then the following new bracket
called the {\em Dirac bracket} on the constrained manifold:
\be
\label{Dirac-def}
\Bigl.\{X_i,X_j\}_D\Bigr|_{C_k=0}:=\{X_i,X_j\}-\{X_i,C_k\}D^{-1}_{k,l}(X)\{X_j,C_l\},
\ee
where, as usual, we imply summations over repeated indices and $D^{-1}_{k,l}$ denotes the $k,l$ entry of $D^{-1}$,  defines a Poisson structure on the constraint manifold.

One, truly remarkable, property of the Dirac brackets is that, provided the matrix $D_{k,l}$ be nondegenerate, the
Dirac bracket (\ref{Dirac-def}) does not depend on an actual parameterization of the constraint submanifold
$\{ C_k=0,\ k=1,\dots,2N\}$: taking instead of $C_k$ {\em arbitrary} functions $f_k(\{C_i\})$ such
that $f_k(\{0\})=0$ and $\Bigl.\det \left[\partial f_k/\partial C_k\right]\Bigr|_{C_\cdot=0}\ne 0$ we obtain the
{\em same} Dirac bracket on the constraint manifold.

We now apply the Dirac bracket procedure to the b.u.t. case. In this case, the constraints are that both $\mathbb A$
and $B\mathbb A {B^{\text{T}}}$ are of the b.u.t. form:
\bea
\label{b.u.t.-constraints-A}
&\mathbb A_{I,J}=0\ \hbox{for}\ I>J, & \det{\mathbb A}_{I,I}=1\\
\label{b.u.t.-constraints-BABt}
&(B\mathbb A {B^{\text{T}}})_{I,J}=0\ \hbox{for}\ I>J, & \det{B\mathbb A {B^{\text{T}}}}_{I,I}=1.
\eea

Observe, first, that since the brackets both between the elements of $\mathbb A$ and
between the elements of $B\mathbb A {B^{\text{T}}}$ are given by the algebra (\ref{Poisson}) and the above
b.u.t. form is a Poisson reduction of this algebra, the brackets between $\mathbb A$-constraints
vanish on the constraint surface (\ref{b.u.t.-constraints-A})
and the same is true for the brackets between $(B\mathbb A {B^{\text{T}}})$-constraints on the constraint
surface (\ref{b.u.t.-constraints-BABt}).
So, the only nonzero brackets on the constraint surface can be an inter-brackets between $\mathbb A$- and
$(B\mathbb A {B^{\text{T}}})$-constraints. This is the point at which the bracket (i) differs from (ii) and (iii).

In the case of bracket (i), because $\mathbb A$ and $B$ commute and because the restriction
(\ref{b.u.t.-constraints-A}) is Poisson for $\mathbb A$-algebra, the inter-brackets between $\mathbb A$- and
$(B\mathbb A {B^{\text{T}}})$-constraints vanish as well, so all the above constraints Poisson commute on the total
constraint surface, but because they do not commute with all other variables, these constraints are not Poisson and {\em the Dirac procedure fails}.\footnote{In principle, again following Dirac's ideology, one might introduce secondary constraints, but this apparently leads to further reduction of the dimension of the
``actual'' phase space, which seems to be not feasible.}

\subsection{Dirac procedure for the upper-triangular case}

We begin with calculating the constraint matrix. We have two sets of constraints,
$$
C_{k,l}=0\ \hbox{and} \ C^\ast_{k,l}=0,\ \hbox{for}\  N\ge k\ge l\ge 0,
$$
where
\be
\begin{array}{lll}
C_{k,l}= & [{\mathbb A}]_{k,l}, \ \hbox{for}\  N\ge k>l\ge 1, & [{\mathbb A}]_{k,k}-1, \ \hbox{for}\ k=1,\dots,N;\\
C^\ast_{k,l}= & \bigl[B{\mathbb A}B^{\text{T}}\bigr]_{k,l} \ \hbox{for}\ N\ge k>l\ge 1,
& \bigl[B{\mathbb A}B^{\text{T}}\bigr]_{k,k}-1, \ \hbox{for}\  k=1,\dots,N;
\end{array}
\label{constr}
\ee
the brackets $\{C,C\}$ and $\{C^\ast,C^\ast\}$ vanish on the constraint surface $C=C^\ast=0$, whereas for
the bracket $\{C,C^\ast\}$ after tedious calculations we obtain the following result (in which we have repeatedly
used already the constraint conditions (\ref{constr})):
\be
\{C_{k,l},C^\ast_{i,j}\}=[B]_{i,k}[B{\mathbb A}^{\text{T}}{\mathbb A}]_{j,l}
+[B{\mathbb A}{\mathbb A}]_{i,l}[B]_{j,k}.
\label{D-matrix}
\ee
We formulate the condition of the nondegeneracy of this matrix in terms of the corresponding system of linear equations: the constraint matrix (\ref{D-matrix}) is nondegenerate iff the matrix equation
\be
P_{-,1}\bigl[B{\mathbb A}{\mathbb A}FB^{\text{T}}+BF^{\text{T}}{\mathbb A}^{\text{T}}{\mathbb A}B^{\text{T}}\bigr]=0
\label{matrix-eq1}
\ee
where $F$ is a nonstrictly upper-triangular matrix, admits only trivial solutions w.r.t. $F$. We can immediately
simplify this equation observing that, for a general upper-triangular $F$, the matrix ${\mathbb A}F$ is also upper
triangular, so instead of solving Eq. (\ref{matrix-eq1}) we can equivalently solve the equation
\be
P_{-,1}\bigl[B{\mathbb A}FB^{\text{T}}+BF^{\text{T}}{\mathbb A}B^{\text{T}}\bigr]=0
\label{matrix-eq2}
\ee
again for an upper-triangular $F$. To simplify this system further, let us introduce the matrix
$g=BF^{\text{T}}B^{-1}$. Then, Eq. (\ref{matrix-eq2}) takes the form
$$
P_{-,1}\bigl[B{\mathbb A}B^{\text{T}}g^{\text{T}}+gB{\mathbb A}B^{\text{T}}\bigr]=0,
$$
which was solved by Bondal in \cite{Bondal}:
the general solution itself can be written in terms of projection operators (see \cite{ChM}):
\be
g=P_{-,1/2}\bigl(B{\mathbb A}B^{\text{T}}\omega_-\bigr)
-P_{+,1/2}\bigl(B{\mathbb A}^{\text{T}}B^{\text{T}}\omega_-^{\text{T}}\bigr),
\ee
where $\omega_-$ is now {\em strictly} lower-triangular matrix. Making the substitution
${\mathbb A}'=B{\mathbb A}B^{\text{T}}$ and $B'=B^{-1}$ and using that
$F^{\text{T}}=B'g{B'}^{-1}$ must be lower triangular,
we can formulate the condition of nondegeneracy
in the following form: given an admissible pair $({\mathbb A}',B')$, i.e., such that $B'{\mathbb A}'{B'}^{\text{T}}$
belongs to the set $\mathcal A$, the equation
\be
P_+[B'g(B')^{-1}]=P_+\Bigl[B'\bigl[P_{-,1/2}({\mathbb A}'\omega_-)
-P_{+,1/2}({{\mathbb A}'}^{\text{T}}\omega_-^{\text{T}})\bigr](B')^{-1}\Bigr]=0
\label{matrix-eq3}
\ee
must have only trivial solutions in the space of strictly lower triangular matrices $\omega_-$.

We first calculate the determinant of the system (\ref{matrix-eq3}) for $B'={\mathbb E}$. In this case, obviously,
this system reads $-P_+({{\mathbb A}}^{\text{T}}\omega_-^{\text{T}})=0$, and the determinant is equal
to the unity irrespectively on ${\mathbb A}$.  This means that for $B$ close to the unit matrix,  this determinant will be nonzero and we can therefore define the Dirac bracket in a layer over the base space $\mathcal A$.  This result establishes a link between this approach and Bondal's one: the neighbourhood of the identity defines the groupoid of morphisms which preserve the space of upper-triangular bilinear forms with $1$ on the diagonal. The pair $(BAB^T, g)$ where $g=BF^{\text{T}}B^{-1}$, belongs to the corresponding Lie algebroid.

We now address the problem of choosing a ``convenient'' parameterization of the $({\mathbb A},B)$-pairs.
In the upper-triangular case, we can express all $a_{i,j}$ with $i<j$ through entries of the matrix $B$. For this,
we write the set of conditions implying that the matrix $B{\mathbb A}B^{\text{T}}$ is upper-triangular:
\be
\sum_{1\le i<j\le n}b_{k,i}b_{l,j}a_{i,j}=-\sum_{s=1}^n b_{k,s}b_{l,s},\quad n\ge k>l\ge 1.
\label{AvisB}
\ee
Additionally, we have a restriction due to Bondal~\cite{Bondal} on the minors $M^{\pm}_{d}$ (see
formulas (\ref{Minor+}) and (\ref{Minor-}) of the
matrix $B$: provided $\det B=1$, $M^{+}_{d}=(-1)^{d(n-d)}M^{-}_{n-d}$, and we assume that all these minors are nonzero.
We have the following technical lemma concerning the determinant of the system (\ref{AvisB}) of linear
equations w.r.t. the entries the matrix ${\mathbb A}$.
\begin{lm}\label{lem-AvisB}
The determinant of the $[n(n-1)/2]\times[n(n-1)/2]$-matrix ${\mathfrak F}_{i<j}^{k>l}:=b_{k,i}b_{l,j}$
is equal to $\prod_{d=1}^{n-1}\bigl[M^{+}_{d}M^{-}_{d}\bigr]$ and is therefore nonzero provided
all upper-right and lower-left minors of the matrix $B$ are nonzero.
\end{lm}
The {\bf proof} uses the Bondal's technique of skew-symmetric forms and can be performed by induction in the
size of the matrix $B$.

We can therefore always express $a_{i,j}$ in terms of $B$ writing
\be
a_{i,j}=F_{i,j}[B].
\label{AF}
\ee
We can also write entries of the transformed matrix $B{\mathbb A}B^{\text{T}}$ in the form
\be
(B{\mathbb A}B^{\text{T}})_{i,j}=BF[B]B^{\text{T}}=\tilde F[B].
\label{BABF}
\ee
Note that the thus defined matrix $\tilde F[B]$ is automatically upper-triangular.

We can therefore replace the set of original constraints (\ref{constr}) by the equivalent set
\be
\begin{array}{ll}
a_{i,j}-F_{i,j}[B]=0, \ i<j, & M^{-}_{d}=(-1)^{d(n-d)}M^{+}_{n-d},\\
a_{k,l}=0,\ k>l, & a_{k,k}=1,\ k=1,\dots,n.
\end{array}
\ee

We now evaluate the brackets between functions
$F_{i,j}[B]$. We begin with case (i) in which the brackets between entries $b_{i,j}$ are
just the Lie--Poisson brackets (\ref{bb-bracket}). In this case,
all constraint equations Poisson commute and ${\mathbb A}$ commutes with $B$, so we have
\be
\{a_{i,j}-F_{i,j}[B],a_{s,p}-F_{s,p}[B]\}=\{a_{i,j},a_{s,p}\}+\{F_{i,j}[B],F_{s,p}[B]\}=0,
\ee
or
\be
\{F_{i,j}[B],F_{s,p}[B]\}=-\{a_{i,j},a_{s,p}\}|_{{\mathbb A}=F[B]}
\ee
(note the minus sign in this relation). The brackets between entries of $F$ are here induced by the
standard Lie--Poisson bracket (\ref{BB-bracket}).

In the case (ii) or (iii), we use the Dirac procedure that implies that brackets
between constraints as well as brackets between constraints and all variables
vanish on the constraint surface. For instance,
$$
\{a_{i,j},a_{s,p}-F_{s,p}[B]\}_D=0,\ \hbox{i.e.,}\ \{a_{i,j},a_{s,p}\}_D=\{a_{i,j},F_{s,p}[B]\}_D,
$$
or, since the $\{a,a\}$ brackets are not changed by the Dirac reduction, we obtain that
\be
\{a_{i,j},F_{s,p}[B]\}_D|_{{\mathbb A}=F[B]}=\{a_{i,j},a_{s,p}\}|_{{\mathbb A}=F[B]}.
\ee
Then for the brackets between entries of the matrix $F$ we obtain
\bea
\{a_{i,j}-F_{i,j}[B],a_{s,p}-F_{s,p}[B]\}_D&=&
\{a_{i,j},a_{s,p}\}|_{{\mathbb A}=F[B]}\nonumber\\
&{}&\quad -2\{a_{i,j},a_{s,p}\}|_{{\mathbb A}=F[B]}
+\{F_{i,j}[B],F_{s,p}[B]\}_D
\nonumber
\eea
and therefore
\be
\{F_{i,j}[B],F_{s,p}[B]\}_D=\{a_{i,j},a_{s,p}\}|_{{\mathbb A}=F[B]}\ \hbox{for}\ i<j\ \hbox{and}\ s<p
\ee
with the plus sign.

{We can also introduce an analogous representation for the matrix $B{\mathbb A}B^{\text{T}}$:
\be
[B{\mathbb A}B^{\text{T}}]_{i,j}:=\left\{\tilde F[B]_{i,j}, i<j; 1, \ i=j; 0, \ i>j\right\}:=\tilde F[B]
\label{tildeF}
\ee
Note that if we impose the standard Lie--Poisson brackets on $B$ then (see the proof and discussion in
Sec.~\ref{s:groupoid})
we obtain that Poisson relations between $\tilde F$ are the same as for ${\mathbb A}$ (with the plus sign)
and $F$ Poisson commute with $\tilde F$.}

{If we again consider the case (ii) or (iii) and apply the Dirac procedure, then, since
the above reasonings remain valid for the set of reparameterized
constraints $[B{\mathbb A}B^{\text{T}}]_{i,j}-\tilde F[B]_{i,j}=0$ and the Dirac brackets between entries
of $B{\mathbb A}B^{\text{T}}$ coincide with the initial brackets,
we obtain that the brackets between entries of $\tilde F$ are
\be
\{\tilde F_{i,j}[B],\tilde F_{s,p}[B]\}_D=\{a_{i,j},a_{s,p}\}|_{{\mathbb A}=\tilde F[B]}.
\ee
}

We can take the set $\{F_{i,j},\tilde F_{s,p}\}$ as {\em new dynamical variables} describing our system;
these variables are not however algebraically independent as they share $[n/2]$ Casimir functions $Y_p$ generated by
$\det(A+\lambda A^{\text{T}})$. But, say, in the case (ii) we also have $[n/2]$ additional Casimir
functions $X_p$
(see Theorem~\ref{th-central-A-B}). Note that for an upper-triangular matrix $\mathbb A$
the principal minors of the matrix ${\mathbb A}^{\text{T}}B^{-1}$ coincide with those of
the matrix $B^{-1}$, so all $X_p$ can be expressed as ratios of principal and upper-right minors of $B$
and are algebraically independent. We can therefore parameterize
the general Poisson leaf of the Dirac Poisson algebra by the variables $F_{i,j}$ and $\tilde F_{s,p}$;
inside each of these two sets the brackets are given by those of the entries of ${\mathbb A}$, the values
of the Casimir functions $Y_p$ coincide for these two sets, and the brackets between $F$ and $\tilde F$ are
determined by the Dirac procedure; we do not evaluate these brackets in this paper and only mention that they are
nontrivial.

\section{Quantization}\label{s:quantization}
\setcounter{equation}{0}

In this section, we quantize the new algebras (\ref{bracket-Q}) for nontrivial $Q$. As above, we concentrate mostly
on the case (ii).

We use the standard trigonometric quantum $R$-matrix
\be
\label{R-q}
R_{12}(q)=\sheet1{\mathbb E}\otimes \sheet2{\mathbb E}
+\sum_{k,l}\sheet1{E}_{k,l}\otimes \sheet2{E}_{l,k}\Bigl[(q-q^{-1})\theta(l-k)+\frac {(q^{1/2}-q^{-1/2})^2}{2}
\delta_{k,l}\Bigr].
\ee
The matrix $R_{12}(q)$ manifests the following useful properties:
\bea
&{}&\bigl[R_{12}(q)\bigr]^{-1}=R_{12}(q^{-1})
\label{R-inverse}\\
&{}&\bigl[R_{12}(q), R^{t_1}_{12}(q)\bigr]=\bigl[R_{12}(q), R^{t_2}_{12}(q)\bigr]=0,
\label{R-comm}\\
&{}&R_{12}(q)+R_{21}(q^{-1})=(q-q^{-1})P_{12}
\label{R-perm}\\
&{}&R_{12}(q)R_{13}(q)R_{23}(q)=R_{23}(q)R_{13}(q)R_{12}(q)\quad\hbox{the Yang--Baxter relation}.
\label{R-YB}
\eea
We use the standard notation: the entries of the matrices ${\mathbb A}=E_{i,j}\otimes a_{i,j}$ and
$B=E_{i,j}\otimes b_{i,j}$ are now operators in the
quantum space. The orders of multiplication in the classical space and in the quantum space can be in principle
different. However, when not stated explicitly, we assume that the order of multiplication in the quantum space
is {\em natural}, i.e., it coincides with the ordering of the matrices ${\mathbb A}$ and $B$ in matrix products.

\begin{lm}\label{lm:quantum}
The quantum version of the case (ii) Poisson algebra reads:
\bea
&{}&
R_{12}(q){\sheet1{\mathbb A}}R_{12}^{t_1}(q){\sheet2{\mathbb A}}
={\sheet2{\mathbb A}}R_{12}^{t_1}(q){\sheet1{\mathbb A}}R_{12}(q)
\label{R-AA}
\\
&{}&
R_{12}(q){\sheet1{B}}{\sheet2{B}}={\sheet2{B}}{\sheet1{B}}R_{12}(q)
\label{R-BB}
\\
&{}&
{\sheet2{\mathbb A}}{\sheet1{B}}R_{12}(q)={\sheet1{B}}R_{12}^{t_2}(q^{-1}){\sheet2{\mathbb A}}
\label{R-AB}
\eea
These commutation relations satisfy the Jacobi relations and ensure the quantum automorphism: the products
$B{\mathbb A}{B^{\text{T}}}$ satisfy the quantum algebra (\ref{R-AA}).
\end{lm}

The {\em proof} is a straightforward but lengthy calculation alongside which we encounter
another Yang--Baxter relation,
\be
R_{23}(q)R_{13}^{t_1}(q^{-1})R_{12}^{t_1}(q^{-1})=R_{12}^{t_1}(q^{-1})R_{13}^{t_1}(q^{-1})R_{23}(q),
\label{YB-new}
\ee
which can be derived from the original relation (\ref{R-YB}) by total transposition in the first space (note that
under this operation entries with $R_{13}$ and $R_{12}$ permute and the entry with $R_{23}$ retains its position).
After that, sandwiching the both sides of the obtained relation between two insertions of $R_{12}(q^{-1})$,
we obtain the original Yang--Baxter relation with the global replacement $q\to q^{-1}$.

Note that it is safe to transpose commutation relations provided we {\em do not} change the order of
multiplication in the quantum space (and entries in the quantum space commute with all $R$-matrix entries).
For instance, from (\ref{R-BB}) we have that
\be
{\sheet1{B^{\text{T}}}}R_{12}^{t_1}(q){\sheet2{B}}={\sheet2{B}}R^{t_1}_{12}(q){\sheet1{B^{\text{T}}}}
\label{R-BB-tr}
\ee
and
\be
{\sheet1{B^{\text{T}}}}{\sheet2{B^{\text{T}}}}R_{12}^{t_1t_2}(q)=R^{t_1t_2}_{12}(q){\sheet2{B^{\text{T}}}}{\sheet1{B^{\text{T}}}}
\label{R-BB-tr-tr}
\ee
and since $R_{12}(q)+R_{21}(q^{-1})=P_{12}$ we can first replace $R_{12}^{t_1t_2}(q)=R_{21}(q)$ by
$-R_{12}(q^{-1})$ and then, multiplying the relation (\ref{R-BB-tr-tr}) by $R_{12}(q)$ from both sides,
we obtain that
\be
R_{12}(q){\sheet1{B^{\text{T}}}}{\sheet2{B^{\text{T}}}}={\sheet2{B^{\text{T}}}}{\sheet1{B^{\text{T}}}}R_{12}(q).
\label{R-BB-tt}
\ee

From (\ref{R-AB}) we have that
\be
{\sheet2{\mathbb A}}R_{12}^{t_1}(q){\sheet1{B^{\text{T}}}}=R_{12}^{t_1t_2}(q^{-1}){\sheet1{B^{\text{T}}}}{\sheet2{\mathbb A}}.
\label{R-AB-tr}
\ee

The proof of the second statement of the lemma (the automorphism) follows from the following chain of equalities:
\bea
&{}&R_{12}(q)\sheet1{B}\sheet1{\mathbb A}\bigl[\sheet1{B^{\text{T}}}R_{12}^{t_1}(q) \sheet2{B}\bigr]
\sheet2{\mathbb A}\sheet2{B^{\text{T}}}
=R_{12}(q)\sheet1{B}\sheet1{\mathbb A}\sheet2{B}R_{12}^{t_1}(q) \sheet1{B^{\text{T}}}\sheet2{\mathbb A}\sheet2{B^{\text{T}}}
\nonumber\\
&{}&
=R_{12}(q)\sheet1{B}\bigl[\sheet1{\mathbb A}\sheet2{B}R_{21}(q)\bigr]
R_{12}^{t_1}(q)\bigl[ R_{21}(q^{-1}) \sheet1{B^{\text{T}}}\sheet2{\mathbb A}\bigr]\sheet2{B^{\text{T}}}
\nonumber\\
&{}&
=\bigl[R_{12}(q)\sheet1{B}\sheet2{B}\bigr]R_{12}^{t_2}(q^{-1})\sheet1{\mathbb A}
R_{12}^{t_1}(q)\sheet2{\mathbb A} R_{12}^{t_1}(q) \sheet1{B^{\text{T}}}\sheet2{B^{\text{T}}}
\nonumber\\
&{}&
=\sheet2{B}\sheet1{B}\bigl[ R_{12}(q) R_{12}^{t_2}(q^{-1})\bigr]\sheet1{\mathbb A}
R_{12}^{t_1}(q)\sheet2{\mathbb A} R_{12}^{t_1}(q) \sheet1{B^{\text{T}}}\sheet2{B^{\text{T}}}
\nonumber\\
&{}&
=\sheet2{B}\sheet1{B} R_{12}^{t_2}(q^{-1})\bigl[ R_{12}(q)\sheet1{\mathbb A}
R_{12}^{t_1}(q)\sheet2{\mathbb A}\bigr] R_{12}^{t_1}(q) \sheet1{B^{\text{T}}}\sheet2{B^{\text{T}}}
\nonumber\\
&{}&
=\sheet2{B}\sheet1{B} R_{12}^{t_2}(q^{-1})\sheet2{\mathbb A}
R_{12}^{t_1}(q)\sheet1{\mathbb A}\bigl[ R_{12}(q) R_{12}^{t_1}(q)\bigr] \sheet1{B^{\text{T}}}\sheet2{B^{\text{T}}}
\nonumber\\
&{}&
=\sheet2{B}\sheet1{B} R_{12}^{t_2}(q^{-1})\sheet2{\mathbb A}
R_{12}^{t_1}(q)\sheet1{\mathbb A}R_{12}^{t_1}(q)\bigl[ R_{12}(q) \sheet1{B^{\text{T}}}\sheet2{B^{\text{T}}} \bigr]
\nonumber\\
&{}&
=\sheet2{B}\bigl[\sheet1{B} R_{12}^{t_2}(q^{-1})\sheet2{\mathbb A} \bigr]
R_{12}^{t_1}(q)\bigl[ \sheet1{\mathbb A}R_{12}^{t_1}(q) \sheet2{B^{\text{T}}}\bigr] \sheet1{B^{\text{T}}} R_{12}(q)
\nonumber\\
&{}&
=\sheet2{B}\sheet2{\mathbb A}\sheet1{B}\bigl[ R_{12}(q)
R_{12}^{t_1}(q)R_{12}(q^{-1})\bigr]\sheet2{B^{\text{T}}}\sheet1{\mathbb A} \sheet1{B^{\text{T}}} R_{12}(q)
\nonumber\\
&{}&
=\sheet2{B}\sheet2{\mathbb A}\bigl[\sheet1{B}
R_{12}^{t_1}(q)\sheet2{B^{\text{T}}}\bigr]\sheet1{\mathbb A} \sheet1{B^{\text{T}}} R_{12}(q)
=\sheet2{B}\sheet2{\mathbb A}\sheet2{B^{\text{T}}}
R_{12}^{t_1}(q)\sheet1{B}\sheet1{\mathbb A} \sheet1{B^{\text{T}}} R_{12}(q).
\nonumber
\eea

We prove the Jacobi relations for the $\{\mathbb A, \mathbb A, B\}$-commutator, the proof for the
$\{\mathbb A, B, B\}$-commutator is analogous.
Making the cyclic change of indices $(1,2,3)\to (2,3,1)$ and doing the total transposition of
(\ref{YB-new}) we obtain that
\be
R_{23}^{t_3}(q^{-1})R_{12}^{t_2}(q^{-1})R_{13}(q)=
R_{13}(q)R_{12}^{t_2}(q^{-1})R_{23}^{t_3}(q^{-1}).
\label{RR-int}
\ee
Using that $R_{23}^{t_3}(q^{-1})+R_{32}^{t_3}(q)=(q-q^{-1})P_{23}^{t_3}$ and that
$$
P_{23}^{t_3}R_{12}^{t_2}(q^{-1})R_{13}(q)=P_{23}^{t_3}R_{13}(q^{-1})R_{13}(q)=P_{23}^{t_3}
$$
and
$$
R_{13}(q)R_{12}^{t_2}(q^{-1})P_{23}^{t_3}=R_{13}(q)R_{13}(q^{-1})P_{23}^{t_3}=P_{23}^{t_3},
$$
we can effectively replace $R_{23}^{t_3}(q^{-1})$ by $R_{32}^{t_3}(q)=R_{23}^{t_2}(q)$ in
(\ref{RR-int}) thus obtaining another Yang--Baxter-type relation
\be
R_{23}^{t_2}(q)R_{12}^{t_2}(q^{-1})R_{13}(q)=
R_{13}(q) R_{12}^{t_2}(q^{-1})R_{23}^{t_2}(q).
\label{YB-another}
\ee
The Jacobi relation now follows from the chain of equalities
\bea
&{}&
\sheet3{\mathbb A}R_{23}^{t_2}(q)\bigl[\sheet2{\mathbb A}\sheet1{B} R_{12}(q)\bigr]R_{13}(q)R_{23}(q)
\nonumber\\
&{}&
=\sheet3{\mathbb A}\sheet1{B}\bigl[ R_{23}^{t_2}(q) R_{12}^{t_2}(q^{-1})R_{13}(q) \bigr]\sheet2{\mathbb A}R_{23}(q)
\nonumber\\
&{}&
=\bigl[ \sheet3{\mathbb A}\sheet1{B}R_{13}(q)\bigr] R_{12}^{t_2}(q^{-1})R_{23}^{t_2}(q)
\sheet2{\mathbb A}R_{23}(q)
\nonumber\\
&{}&
=\bigl[ \sheet3{\mathbb A}\sheet1{B}R_{13}(q)\bigr] R_{12}^{t_2}(q^{-1})R_{23}^{t_2}(q)
\sheet2{\mathbb A}R_{23}(q)
\nonumber\\
&{}&
=\sheet1{B}R_{13}^{t_3}(q^{-1}) R_{12}^{t_2}(q^{-1})\bigl[\sheet3{\mathbb A}R_{23}^{t_2}(q)
\sheet2{\mathbb A}R_{23}(q)\bigr]
\nonumber\\
&{}&
=\sheet1{B}\bigl[ R_{13}^{t_3}(q^{-1}) R_{12}^{t_2}(q^{-1})
R_{23}(q)\bigr]\sheet2{\mathbb A} R_{23}^{t_2}(q)\sheet3{\mathbb A}
\nonumber\\
&{}&
=R_{23}(q)\bigl[  \sheet1{B} R_{12}^{t_2}(q^{-1})\sheet2{\mathbb A} \bigr] R_{13}^{t_3}(q^{-1})
R_{23}^{t_2}(q)\sheet3{\mathbb A}
\nonumber\\
&{}&
=R_{23}(q)\sheet2{\mathbb A}  \sheet1{B}\bigl[  R_{12}(q) R_{13}^{t_3}(q^{-1})
R_{23}^{t_2}(q)\bigr]\sheet3{\mathbb A}
\nonumber\\
&{}&
=R_{23}(q)\sheet2{\mathbb A}R_{23}^{t_2}(q)\bigl[ \sheet1{B} R_{13}^{t_3}(q^{-1})
\sheet3{\mathbb A}\bigr] R_{12}(q)
\nonumber\\
&{}&
=\bigl[ R_{23}(q)\sheet2{\mathbb A}R_{23}^{t_2}(q)\sheet3{\mathbb A}\bigr] \sheet1{B} R_{13}(q)
R_{12}(q)
\nonumber\\
&{}&
= \sheet3{\mathbb A}R_{23}^{t_2}(q)\sheet2{\mathbb A} \sheet1{B}R_{23}(q) R_{13}(q)
R_{12}(q),
\nonumber
\eea
and using the standard Yang--Baxter relation (\ref{R-YB}) we obtain the identity. The lemma is proved.

\section{Quantum affine algebras}\label{s:quant-spectral}
\setcounter{equation}{0}
Similarly to the approach of paper \cite{ChM} we can introduce the affine version of the quantum relations.
We now have infinite-dimensional quantum algebras with the generators
\be
{\mathbb A}(\lambda)=\sum_{i,j}\sum_{k=0}^\infty E_{i,j}\otimes a_{i,j}^{(k)}\lambda^{-k},\qquad
B(\lambda)=\sum_{i,j}\sum_{k=0}^\infty E_{i,j}\otimes b_{i,j}^{(k)}\lambda^{-k}.
\ee
We also have the quantum affine $R$-matrix (depending on the spectral parameters $\lambda$ and $\mu$):
\bea
&{}&
R_{12}(\lambda,\mu;q):=\frac{\lambda-\mu}{q^{-1}\lambda-q\mu}\sum_{i\ne j}\sheet1{E}_{i,i}\otimes \sheet2{E}_{j,j}
+\sum_i \sheet1{E}_{i,i}\otimes \sheet2{E}_{i,i}
\nonumber\\
&{}&\qquad
+\frac{(q^{-1}-q)\lambda}{q^{-1}\lambda-q\mu}\sum_{i<j}\sheet1{E}_{i,j}\otimes \sheet2{E}_{j,i}
+\frac{(q^{-1}-q)\mu}{q^{-1}\lambda-q\mu}\sum_{i>j}\sheet1{E}_{i,j}\otimes \sheet2{E}_{j,i}.
\label{R-MN}
\eea
We choose the normalization such that
\be
\bigl[R_{12}(\lambda,\mu;q)\bigr]^{-1}=R_{12}(\lambda,\mu;q^{-1}).
\label{R-MN-inverse}
\ee
We also have several useful relations:
\bea
&{}&
R_{12}(\lambda,\mu;q)=R_{12}^{t_1t_2}(\lambda^{-1},\mu^{-1};q^{-1}),\label{R-1}\\
&{}&
R_{12}(\lambda,\mu;q)=R_{12}(\mu^{-1},\lambda^{-1};q),\label{R-2}\\
&{}&
R_{12}(\lambda,\mu;q)R_{12}^{t_1}(\rho,\nu; q^{-1}) R_{12}(\lambda,\mu;q^{-1})=R_{12}^{t_1}(\rho,\nu; q^{-1}).
\label{R-3}
\eea

The $R$-matrix $R_{12}(\lambda,\mu;q)$ satisfies the Yang--Baxter relations
\be
R_{12}(\lambda,\mu;q)R_{13}(\lambda,\rho;q)R_{23}(\mu,\rho;q)=
R_{23}(\mu,\rho;q)R_{13}(\lambda,\rho;q)R_{12}(\lambda,\mu;q)
\label{YB-MN}
\ee

The main lemma reads

\begin{lm}\label{lm:quantum-mn}
The affine quantum version of the case (ii) Poisson algebra reads:
\bea
&{}&
R_{12}(\lambda,\mu;q){\sheet1{\mathbb A}}(\lambda)R_{12}^{t_1}(\lambda^{-1},\mu;q){\sheet2{\mathbb A}}(\mu)
={\sheet2{\mathbb A}}(\mu)R_{12}^{t_1}(\lambda^{-1},\mu;q){\sheet1{\mathbb A}}(\lambda)R_{12}(\lambda,\mu;q)
\label{R-AA-mn}
\\
&{}&
R_{12}(\lambda,\mu;q){\sheet1{B}}(\lambda){\sheet2{B}}(\mu)
={\sheet2{B}}(\mu){\sheet1{B}}(\lambda)R_{12}(\lambda,\mu;q)
\label{R-BB-mn}
\\
&{}&
{\sheet2{\mathbb A}}(\mu){\sheet1{B}}(\lambda)R_{12}(\lambda,\mu;q)
={\sheet1{B}}(\lambda)R_{12}^{t_2}(\mu,\lambda^{-1};q^{-1}){\sheet2{\mathbb A}}(\mu)
\label{R-AB-mn}
\eea
Note that in the last relation we can equivalently substitute
$$
R_{12}^{t_2}(\mu,\lambda^{-1};q^{-1})=R_{12}^{t_2}(\lambda,\mu^{-1};q^{-1}).
$$
These commutation relations satisfy the Jacobi relations and ensure the quantum automorphism: the products
$$
B(\lambda){\mathbb A}(\lambda){B^{\text{T}}}(\lambda^{-1})
$$
satisfy the quantum algebra (\ref{R-AA-mn}).
\end{lm}

The {\em proof} is straightforward. We have
a Yang--Baxter-type relation
\bea
&{}&
R_{23}(\mu,\nu;q)R_{12}^{t_2}(\lambda,\mu^{-1};q^{-1})R_{13}^{t_3}(\lambda,\nu^{-1};q^{-1})
\nonumber\\
&{}&\qquad\qquad =R_{13}^{t_3}(\lambda,\nu^{-1};q^{-1})R_{12}^{t_2}(\lambda,\mu^{-1};q^{-1})R_{23}(\mu,\nu;q),
\label{YB-new-mn}
\eea
which can again be derived from the standard relation (\ref{YB-MN}) if we perform the transposition in the first
space (then, again, items with $R_{13}$ and $R_{12}$ permute) and take relation (\ref{R-1}) into account.
After this, the two chains of calculations are completely analogous to those in the proof of Lemma~\ref{lm:quantum}; we just need
to set the spectral parameters in the proper places and use the corresponding identities. For example, we obtain the
commutation relations between $\sheet1{B^{\text{T}}}(\lambda^{-1})$ and $\sheet2{B^{\text{T}}}(\mu^{-1})$  by transposing
relations (\ref{R-BB-mn}) and substituting $(\lambda, \mu)\to (\lambda^{-1}, \mu^{-1})$, which gives
$$
\sheet1{B^{\text{T}}}(\lambda^{-1})\sheet2{B^{\text{T}}}(\mu^{-1})R_{12}^{t_1t_2}(\lambda^{-1},\mu^{-1};q)=
R_{12}^{t_1t_2}(\lambda^{-1},\mu^{-1};q)\sheet2{B^{\text{T}}}(\mu^{-1})\sheet1{B^{\text{T}}}(\lambda^{-1}).
$$
But
$$
R_{12}^{t_1t_2}(\lambda^{-1},\mu^{-1};q)=R_{12}(\lambda,\mu;q^{-1})=\bigl[R_{12}(\lambda,\mu;q)\bigr]^{-1},
$$
so multiplying from the both sides by $R_{12}(\lambda,\mu;q)$, we obtain the familiar relation
$$
R_{12}(\lambda,\mu;q)\sheet1{B^{\text{T}}}(\lambda^{-1})\sheet2{B^{\text{T}}}(\mu^{-1})=
\sheet2{B^{\text{T}}}(\mu^{-1})\sheet1{B^{\text{T}}}(\lambda^{-1})R_{12}(\lambda,\mu;q).
$$
We omit the proof of the algebra automorphism as this proof repeats that in Lemma~\ref{lm:quantum}; we only present
the proof of the Jacobi relations for the $\{\mathbb A(\lambda), B(\mu), B(\nu)\}$-triple relations:
\bea
&{}&
\sheet1{\mathbb A}(\lambda)
\bigl[\sheet3{B}(\nu)\sheet2{B}(\mu)R_{23}(\mu,\nu;q)\bigr]R_{21}(\mu,\lambda;q)R_{31}(\nu,\lambda;q)
\nonumber\\
&{}&
=R_{23}(\mu,\nu;q)\bigl[\sheet1{\mathbb A}(\lambda)\sheet2{B}(\mu)R_{21}(\mu,\lambda;q)\bigr]
\sheet3{B}(\nu)R_{31}(\nu,\lambda;q)
\nonumber\\
&{}&
=R_{23}(\mu,\nu;q)\sheet2{B}(\mu)
R_{12}^{t_2}(\lambda,\mu^{-1};q^{-1})\bigl[\sheet1{\mathbb A}(\lambda)\sheet3{B}(\nu)
R_{31}(\nu,\lambda;q)\bigr]
\nonumber\\
&{}&
=\bigl[R_{23}(\mu,\nu;q)\sheet2{B}(\mu)\sheet3{B}(\nu)\bigr]
R_{12}^{t_2}(\lambda,\mu^{-1};q^{-1})
R_{13}^{t_3}(\lambda,\nu^{-1};q^{-1})\sheet1{\mathbb A}(\lambda)
\nonumber\\
&{}&
=\sheet3{B}(\nu)\sheet2{B}(\mu)
\bigl[R_{23}(\mu,\nu;q)R_{12}^{t_2}(\lambda,\mu^{-1};q^{-1})R_{13}^{t_3}(\lambda,\nu^{-1};q^{-1})\bigr]
\sheet1{\mathbb A}(\lambda)
\nonumber\\
&{}&
=\sheet3{B}(\nu)R_{13}^{t_3}(\lambda,\nu^{-1};q^{-1})\bigl[\sheet2{B}(\mu)
R_{12}^{t_2}(\lambda,\mu^{-1};q^{-1})
\sheet1{\mathbb A}(\lambda)\bigr]R_{23}(\mu,\nu;q)
\nonumber\\
&{}&
=\bigl[\sheet3{B}(\nu)R_{13}^{t_3}(\lambda,\nu^{-1};q^{-1})\sheet1{\mathbb A}(\lambda)\bigr]
\sheet2{B}(\mu)R_{21}(\mu,\lambda;q)R_{23}(\mu,\nu;q)
\nonumber\\
&{}&
=\sheet1{\mathbb A}(\lambda)\sheet3{B}(\nu)\sheet2{B}(\mu)
R_{31}(\nu,\lambda,;q)R_{21}(\mu,\lambda;q)R_{23}(\mu,\nu;q),
\nonumber
\eea
and we again obtain the standard Yang--Baxter relation in the first and last terms. The Jacobi relation is
therefore proved.

\section{Poisson groupoid structure}\label{s:groupoid}
\setcounter{equation}{0}

In this section, we investigate the possibility of inducing brackets on the upper-triangular matrix ${\mathbb A}$
from the Lie--Poisson brackets (\ref{Poisson-b}) on $B$. We shall demonstrate that the obtained
brackets then satisfy the definition of a Poisson (symplectic) groupoid~\cite{Weinstein},~\cite{Mikami-Wein}.

As we have already demonstrated in Sec.~\ref{s:reductions}, we can express entries of ${\mathbb A}$ as the functions
$F_{i,j}[B]$; we have already shown that the brackets between $F_{i,j}[B]$ and $F_{s,p}[B]$
are given by minus the bracket between $a_{i,j}$ and $a_{s,p}$ upon the substitution ${\mathbb A}=F[B]$. We now
evaluate all other brackets of this system.

We begin with evaluating the bracket between $F_{i,j}[B]$ and $b_{s,p}$ induced by the bracket
(\ref{Poisson-b}). Note first that
\bea
\{b_{k,i}b_{l,j},b_{s,p}\}&=&\theta(i-p)[b_{k,p}b_{l,j}]b_{s,i}-\theta(s-k)[b_{s,i}b_{l,j}]b_{k,p}
\nonumber\\
&{}&-\theta(p-j)[b_{k,i}b_{l,p}]b_{a,j}+\theta(l-s)[b_{k,i}b_{s,j}]b_{l,p}.
\label{B1}
\eea
In this expression we group in the square brackets in the r.h.s. the terms that are entries of
the matrix ${\mathfrak F}_{\alpha<\beta}^{\kappa>\lambda}:=b_{\kappa,\alpha}b_{\lambda,\beta}$ provided the original
combintation
$b_{k,i}b_{l,j}$ was an entry of this matrix. If we take $i=j$ in (\ref{B1}) then in the r.h.s. we have either terms
with coinciding right indices or terms that are again entries of ${\mathfrak F}$. We are therefore able
to evaluate explicitly the Poisson bracket of $b_{s,p}$ with
$$
F_{i,j}[B]=\sum_{k>l}\bigl[{\mathfrak F}^{-1}]_{i<j}^{k>l}\Bigl[\sum_s b_{k,s}b_{l,s}\Bigr].
$$
After a tedious algebra we obtain a compact answer, which can be written in a convenient $r$-matrix form:
taking ${\mathbb F}:=\{F_{i,j}[B], \ i<j; \ 1, \ i=j, \ 0, \ i>j\}$ to be an upper-triangular matrix, we
merely obtain that
\be
\{\sheet1 B\ocomma \sheet2 {\mathbb F}\}=\sheet1 B r_{12} \sheet2 {\mathbb F}
-\sheet1 B\sheet2 {\mathbb F} r_{12}^{t_1},
\label{PLB-B-F}
\ee
or, in the component form,
\be
\{b_{i,j},f_{k,l}\}=\sum_s \theta(s-j)\delta_{k,j}b_{i,s}f_{s,l}-\sum_s \theta(j-s)\delta_{j,l}b_{i,s}f_{k,s}.
\label{PLB-B-F-1}
\ee

We can take the above formula as the {\it definition} of a new bracket between
$B$ and ${\mathbb F}$ assuming that ${\mathbb F}$ have the Poisson bracket of the form (\ref{Poisson-r})
{\it with the overall negative sign},
\be
\{\sheet1{\mathbb F}\ocomma\sheet2{\mathbb F}\}=-r_{12}(\sheet1{\mathbb F}\sheet2{\mathbb F})
+(\sheet1{\mathbb F}\sheet2{\mathbb F})r_{12}-\sheet1{\mathbb F}r_{12}^{t_1}\sheet2{\mathbb F}
+\sheet2{\mathbb F}r_{12}^{t_1}\sheet1{\mathbb F},
\label{B2}
\ee
and $B$ have the standard brackets (\ref{Poisson-b}). If we then perform the
mapping ${\mathbb F}\mapsto \tilde {\mathbb F}:=B{\mathbb F}B^{\text{T}}$ we obtain
using only (\ref{Poisson-b}), (\ref{PLB-B-F}), and (\ref{B2})
the following Poisson relations:
\bea
&{}&
\{\sheet1{\tilde{\mathbb F}}\ocomma\sheet2{\tilde{\mathbb F}}\}=
r_{12}(\sheet1{\tilde{\mathbb F}}\sheet2{\tilde{\mathbb F}})
-(\sheet1{\tilde{\mathbb F}}\sheet2{\tilde{\mathbb F}})r_{12}
+\sheet1{\tilde{\mathbb F}}r_{12}^{t_1}\sheet2{\tilde{\mathbb F}}
-\sheet2{\tilde{\mathbb F}}r_{12}^{t_1}\sheet1{\tilde{\mathbb F}},
\label{tilde-A}\\
&{}&\{\sheet1{B}\ocomma\sheet2{\tilde{\mathbb F}}\}=
r_{12}(\sheet1{B}\sheet2{\tilde{\mathbb F}})
-\sheet2{\tilde{\mathbb F}}r_{12}^{t_1}\sheet1{B}.
\label{tilde-B}
\eea
The component form of writing of the latter relation is
\be
\{b_{i,j},\tilde f_{k,l}\}=\theta(i-k)b_{k,j}\tilde f_{i,l}-\theta(l-i)b_{l,j}\tilde f_{k,i}.
\label{tilde-B-1}
\ee

Note first that this new Poisson algebra of ${\mathbb F}$ and $B$ satisfies all the Jacobi relations even if we
consider {\em independent} pairs $({\mathbb F},B)\in GL_N\times GL_N$. Actually, the most
economic way to see the
satisfaction of the Poisson Jacobi relations
is by {\em quantizing} the corresponding algebra; the corresponding
{\em quantum commutation relations} are
\be
\sheet2{\tilde{\mathbb F}}R^{t_1}_{12}[q]\sheet1B=R_{12}[q]\sheet1 B\sheet2{\tilde{\mathbb F}}.
\label{quantum}
\ee

Second, this algebra admits Poisson reduction to any block-upper triangular form: all the constraints
${\mathbb F}_{I,J}=0$ for $I>J$ and $\tilde{\mathbb F}_{I,J}=0$ for $I>J$ are now Poissonnian.
Third, evaluating the
brackets between ${\mathbb F}$  and $\tilde{\mathbb F}$ we merely obtain that
$$
\{\sheet1{\tilde{\mathbb F}}\ocomma\sheet2{{\mathbb F}}\}=0,
$$
so these algebras are in fact totally separated.

These two conditions: that ${\mathbb F}\to \tilde {\mathbb F}=B{\mathbb F} B^{\text{T}}$ is an anti-Poisson
mapping and that ${\mathbb F}$ Poisson commutes with $\tilde {\mathbb F}$ imply that these Poisson brackets
are in fact those for a symplectic groupoid \cite{Weinstein}. We can therefore formulate the following lemma
and conjecture.

\begin{lm}\label{lm:groupoid}
Considering a pair $({\mathbb F},B)\in GL_N\times GL_N$ endowed with the brackets
(\ref{Poisson-b}), (\ref{PLB-B-F}), and (\ref{B2}) we have that
\begin{enumerate}
\item[(i)] the source $s:({\mathbb F},B)\to {\mathbb F}$ and target $t:({\mathbb F},B)\to \tilde{\mathbb F}
\equiv B{\mathbb F} B^{\text{T}}$ mappings are correspondingly antiautomorphism and automorphism of the
Poisson algebra (\ref{Poisson-r});
\item[((ii)] the above algebraic elements ${\mathbb F}$ and $\tilde{\mathbb F}$ Poisson commute;
\item[(iii)] restrictions of ${\mathbb F}$ and $\tilde{\mathbb F}$ to any b.u.t. form are Poissonnian
w.r.t. the total Poisson algebra of the pair $({\mathbb F},B)$;
\item[(iv)] in the upper-triangular case, the Poisson algebra for ${\mathbb F}={\mathbb A}$ is induced by the
Lie--Poisson brackets (\ref{BB-bracket}) for $B$;
\item[(v)] the quantum version of the Poisson $({\mathbb F},B)$-algebra is given by (\ref{R-BB}),
(\ref{quantum}), and the relation
\be
R_{12}(q^{-1}){\sheet1{\mathbb F}}R_{12}^{t_1}(q^{-1}){\sheet2{\mathbb F}}
={\sheet2{\mathbb F}}R_{12}^{t_1}(q^{-1}){\sheet1{\mathbb F}}R_{12}(q^{-1})
\label{R-AA-inverse}
\ee
inverse to (\ref{R-AA}).
\end{enumerate}
\end{lm}

We also observe that the above brackets satisfy the definition of the Poisson (symplectic) groupoid
(see \cite{Weinstein}, \cite{Mikami-Wein})

\begin{defin}\label{def:Poisson}
A {\em Poisson} (symplectic) {\em groupoid} $\Gamma$ is a symplectic manifold $(\Gamma,\Omega)$, where
$\Omega$ is a Poisson structure such that the graph ${\mathcal M}:=\{(x,y,m(x,y))\in 
\Gamma\times\Gamma\times\Gamma| (x,y)\in \Gamma_2\}$ of the groupoid multiplication $m$ is a Lagrangian
submanifold of
$$
(\Gamma,\Omega)\times (\Gamma,\Omega)\times (\Gamma, -\Omega).
$$
\end{defin}

In other words, given three mutually Poisson commuting pairs $({\mathbb F}_i,B_i)\in GL_N\times GL_N$, $i=1,2,3$ 
each endowed with a Poisson structure $\pm\Omega$ with the plus signs for the first two pairs and the minus sign for the
third pair, we obtain the structure of the Poisson groupoid if the following three sets of constraints are Lagrangian:
\be
f:=B_3-B_2B_1,\quad g:={\mathbb F}_3-{\mathbb F}_1,\quad h:={\mathbb F}_2-B_1{\mathbb F}_1B_1^{\text{T}},
\label{fgh}
\ee
i.e., all these constraints Poisson commute on the constraint surface.

\begin{lm}\label{lm:Poisson-groupoid}
The pair $\Gamma:=({\mathbb F},B)$ endowed with the Poisson bracket $\Omega$ defined by 
the Poissom relations (\ref{Poisson-b}), (\ref{PLB-B-F}), and (\ref{B2}) satisfies Definition~\ref{def:Poisson}
of the Poisson (symplectic) groupoid.
\end{lm}

The {\it proof} is the direct calculation: taking three mutually commuting pairs $({\mathbb F}_i,B_i)$ endowed with 
the respective Poisson, Poisson, and anti-Poisson structures, we can easily verify that all six brackets
$\{f,f\}$, $\{g,g\}$, $\{h,h\}$, $\{f,g\}$, $\{f,h\}$, $\{g,h\}$ vanish on the constraint surface.
Note that in order for the last bracket $\{g,h\}=\{{\mathbb F}_1,\tilde{\mathbb F}_1\}$ to vanish, the source and
target projections of the pair $({\mathbb F},B)$ must Poisson commute.

\begin{conjecture}\label{con:symplectic-groupoid}
Lemma~\ref{lm:Poisson-groupoid} implies the existence of the symplectic form
on the pairs $({\mathbb F},B)$ for the b.u.t.-restricted matrices ${\mathbb F}$. This symplectic
structure must generalize the Bondal's symplectic structure \cite{Bondal} on the pairs $({\mathbb A},B)$
in the upper-triangular case.
\end{conjecture}

Note that the mapping ${\mathbb F}\mapsto B{\mathbb F}B^{\text{T}}$
is now an antiautomorphism of the Poisson algebra for ${\mathbb F}$. This complies with the
Poisson groupoid construction but is opposite to the ideology of Poisson symmetric spaces,
which we advocate in the main part of the paper. One reason for us to ``dislike'' an otherwise nice
symplectic groupoid construction is that because of the total separation of variables ${\mathbb F}$ and ${\tilde{\mathbb F}}$,
the dynamics described by the mapping ${\mathbb F}\mapsto B{\mathbb F}B^{\text{T}}$
becomes trivial.\footnote{A minor subtlety appear when concerning central elements of this algebra: ${\mathbb F}$ and ${\tilde{\mathbb F}}$ share the same set of
elements $Y_p$ that are central for both these sets; however the original
Poisson--Lie bracket for $B$ has the full Poisson dimension $n(n-1)$, so we must
have additional $[n/2]$ elements $Q_p$ that are, first,
algebraically independent with ${\mathbb F}$ and ${\tilde{\mathbb F}}$ and, second,
do not commute with $Y_p$. The constraints $Y_p=c_p,\ Q_p=0$, $p=1,\dots, [n/2]$, are then of the second kind, the matrix $\{Y_p,Q_r\}$ is nondegenerate,
and we can again implement the Dirac procedure w.r.t. these constraints; this procedure does not change the
brackets inside the set $\{{\mathbb F}, \tilde{\mathbb F}\}$ of remaining dynamical variables.}

\vskip 2mm \noindent{\bf Acknowledgements.}  The authors are
specially grateful to Rui Fernandes for his insights on Poisson symmetric spaces and other helpful comments
and to Alexei Bondal for his numerous and helpful explanations of the symplectic groupoid construction.
We would also like to thank Alexander Molev, Stefan Kolb, and Kirill Mackenzie for many
enlighting conversations.

The work of L.Ch. was supported in part by
the Russian Foundation for Basic Research (Grant Nos. 14-01-00860-a and 13-01-12405-ofi-m),
by the Program Mathematical Methods for Nonlinear Dynamics, and partly by the  Engineering and Physics Sciences Research Council EP/J007234/1.

The work of M. Mazzocco was supported  by the Engineering and Physics Sciences Research Council EP/J007234/1.



\def\thetheorem{\Alph{section}.\arabic{theorem}}
\def\theprop{\Alph{section}.\arabic{prop}}
\def\thelemma{\Alph{section}.\arabic{lm}}
\def\thecor{\Alph{section}.\arabic{cor}}
\def\theexam{\Alph{section}.\arabic{exam}}
\def\theremark{\Alph{section}.\arabic{remark}}
\def\theequation{\Alph{section}.\arabic{equation}}

\setcounter{section}{0}

\appendix{Notation}\label{se:notation}
\setcounter{equation}{0}

In this paper we use the standard notation
$$
\sheet1{\mathbb A}= \mathbb A\otimes E, \quad\hbox{and}\quad
\sheet2{\mathbb A}=E\otimes\mathbb A,
$$
where $E$ is the $N\times N$ identity matrix. The expression $r_{12}^{t_i}$ denotes the matrix (\ref{r-matrix}) transposed w.r.t. the arguments of the $i$th space ($i=1,2$). The $r$-matrix satisfies the classical Yang--Baxter equation ensuring the
Jacobi identities for the Poisson brackets and the additional property that
\be
r_{12}+r_{12}^{t_1t_2}=2P_{12},
\ee
where $P_{12}=\sum_{i,j}\sheet1{E}_{i,j}\otimes \sheet2{E}_{j,i}$ is the standard permutation $r$-matrix satisfying
\be
P_{12} (\sheet1{A}\otimes \sheet2{B})=(\sheet1{B}\otimes \sheet2{A})P_{12}
\ee
for any matrices $A$ and $B$.

In this notation, $\{\sheet1{\mathbb A}\ocomma\sheet2{\mathbb A}\}$ is a tensor of $4$ components. To extract the bracket between entries $a_{i,j}$ and $a_{k,l}$ of the matrix $\mathbb A$ we need to compute the $^{ik}_{jl}$ component of the tensor $r_{12}(\sheet1{\mathbb A}\otimes\sheet2{\mathbb A})
-(\sheet1{\mathbb A}\otimes\sheet2{\mathbb A})r_{12}+\sheet1{\mathbb A}r_{12}^{t_1}\sheet2{\mathbb A}
-\sheet2{\mathbb A}r_{12}^{t_1}\sheet1{\mathbb A}$. This gives (\ref{Poisson}).


\begin{thebibliography}{99}

\footnotesize\itemsep=0pt

\bibitem{Bondal}
A.~Bondal, {\it A symplectic groupoid of
triangular bilinear forms and the
braid groups}, preprint IHES/M/00/02 (Jan. 2000); {\sl Izv. Math.}, {\bf 68} (2004) 659--708.

\bibitem{Bon1}
A.~Bondal, {\it Symplectic groupoids related to Poisson--Lie groups},
{\it Tr. Mat. Inst. Steklova,}\/ , {\bf 246} (2004) 43--63.

\bibitem{ChM} L. Chekhov, M. Mazzocco, {\it Block triangular bilinear forms and braid group action},
{\it Comm. Math. Phys.}, {\bf 322}, (2013) no.1:49--71. 


\bibitem{ChM1} L. Chekhov, M. Mazzocco, {\it Isomonodromic deformations and twisted Yangians arising in Teichm\"uller theory},
{\sl Adv. Math.},\/{\bf 226(6)} (2011) 4731-4775.

\bibitem{CF1}
M.~Crainic and R.~Fernandes,
Integrability of Lie brackets.
{\it Ann. of Math.}\/ {\bf 157} (2003), no. 2, 575--620.

\bibitem{Dirac}
P.~A.~M. Dirac, Generalized Hamiltonian dynamics, {\it Canadian J. Math.} {\bf2} (1950) 129-148.

\bibitem{Dr}
V.~G. Drinfeld,
On Poisson Homogeneous Spaces of Poisson--Lie Groups,
{\it Theoret. and Math. Phys.,}\/ {\bf 95}, (1993), 524--525.


\bibitem{F2}
Fernandes R.~L., A note on Poisson Symmetric Spaces, {\it Proceedings of the Cornelius Lanczos International Centenary Conference},\/ Eds. J. Davis Brown, Moody T. Chu, Donald C. Ellison, Robert J. Plemmons, SIAM Philadelphia, USA (1994) 638--642.

\bibitem{FM}
V.~Fock and A.~Marshakov,
A note on quantum groups and relativistic Toda theory,
{\it Nucl. Phys. B.}\/ {\bf 56} (1997) 208--214.

\bibitem{GK}
A.~M.~Gavrilik and A.~U.~Klimyk, $q$-Deformed orthogonal and pseudo-orthogonal algebras and their representations,
{\it Lett. Math. Phys.} {\bf21} (1991) 215--220. 

\bibitem{GY}
K.~R. Goodearl and M.~Yakimov, Poisson structures on affine spaces and flag varieties. II. {\it Trans. Amer. Math. Soc.}\/ {\bf 361} (2009), no. 11:5753--5780.



\bibitem{Karasev}
M.~V.~karasev, Analogues of objects of Lie group theory by nonlinear Poisson brackets, {\it Math. USSR Izvestia}
{\bf 28} (1987) 497--527.

\bibitem{kirill}
Mackenzie, Kirill,
General Theory of Lie Groupoids and Lie Algebroids,
{\it LMS Lect. Note Series}\/ {\bf 213} (2005).

\bibitem{Lu}
Lu J.~H., Classical Dynamical $r$-Matrices
and   Poisson Structures on $G\slash H$ and $K\slash T$,
{\it Commun. Math. Phys.,}\/ {\bf 212}, (2000) 337--370.

\bibitem{LuW}
Lu J.~H. and Weinstein A., Poisson Lie groups, dressing transformations and Bruhat decompositions,
{\it J. Diff. Geom.,} {\bf 31} (1990) 501--526.


\bibitem{Mikami-Wein}
K. Mikami and A. Weinstein, Moments and reductions for symplectic groupoids, {\sl Publ. RIMS, Kyoto Univ.}
{\bf 24} (1988) 121--140.

\bibitem{Molev}
Molev A., Yangians and classical Lie algebras.
{\it Mathematical Surveys and Monographs,}\/ {\bf 143}, American Mathematical Society, Providence, RI, (2007).


\bibitem{MR} A. Molev, E. Ragoucy, {\it Symmetries and invariants of twisted quantum algebras and associated Poisson algebras},
{\sl Rev. Math. Phys.}, {\bf 20}(2) (2008) 173--198.

\bibitem{MRS} A. Molev, E. Ragoucy, P. Sorba, {\it Coideal subalgebras in quantum affine algebras}, {\sl Rev. Math. Phys.}, {\bf 15} (2003) 789--822.


\bibitem{NR} 
Nelson J.E., Regge T., Homotopy groups and $(2{+}1)$-dimensional
quantum gravity, {\it Nucl. Phys.~B} {\bf 328} (1989), 190--199.

\bibitem{NRZ}
Nelson J.E., Regge T., Zertuche F., Homotopy groups and
$(2+1)$-dimensional quantum de~Sitter gravity, {\it Nucl. Phys.~B} {\bf 339}
(1990), 516--532.





\bibitem{Wein}
S. Weinberg, The Quantum Theory of Fields, Vol.~1, Cambridge Univ. Press, 1995.

\bibitem{Weinstein}
A. Weinstein, Coisotropic calculus and Poisson groupoids, {\it J. Math. Soc. Japan} {\bf 40}(4) (1988) 705--727.

\end{thebibliography}
\end{document}